\documentclass[12pt]{amsart}
\usepackage{amsfonts,amssymb,latexsym,amsmath,amsthm,color,dsfont}
\usepackage{enumerate,enumitem}
\usepackage{fullpage}
\usepackage{multirow}
\usepackage{makecell,array}

\usepackage{graphicx,anysize,float,booktabs,geometry}
\usepackage{caption} 
\usepackage{hyperref}

\hypersetup{citecolor=red, linkcolor=blue, colorlinks=true}

\marginsize{2.7cm}{2.7cm}{2cm}{0cm}
\allowdisplaybreaks

\newcommand\F{\mathbb{F}}

\newcommand{\mc}{\mathcal}

\newcommand{\cL}{\mathcal{L}}

\newcommand{\cP}{\mathcal{P}}

\newcommand\cS{\mathcal{S}}

\newcommand\A{{\rm A}}
\newcommand\Aut{{\rm Aut}}
\newcommand\C{{\rm C}}
\newcommand\diag{{\sf diag}}
\newcommand\antidiag{{\sf antidiag}}
\newcommand\sym{{\rm S}}

\newcommand\rQ{{\rm Q}}

\newcommand{\Fix}{{\rm Fix}}
\newcommand{\GQ}{{\rm GQ}}

\newcommand\PGU{{\rm PGU}}
\newcommand\PSU{{\rm PSU}}

\newcommand\ASL{{\rm ASL}}
\newcommand\PSL{{\rm PSL}}
\newcommand\GU{{\rm GU}}
\newcommand\SU{{\rm SU}}
\newcommand\SO{{\rm SO}}
\newcommand\SL{{\rm SL}}
\newcommand\GL{{\rm GL}}

\newcommand\w{{\omega}}
\newcommand\ld{{\lambda}}

\theoremstyle{plain}
\newtheorem{theorem}{Theorem}[section]
\newtheorem{problem}[theorem]{Problem}
\newtheorem{lemma}[theorem]{Lemma}
\newtheorem{corollary}[theorem]{Corollary}

\numberwithin{equation}{section}

\theoremstyle{remark}

\def\<{\langle}
\def\>{\rangle}

\def\w{\omega}
\def\a{\alpha}
\def\b{\beta}

\title{Nonexistence of generalized quadrangles admitting a point-primitive and line-primitive automorphism group with socle $\PSU(3,q)$, $q\geq 3$}

\author[Lu, Zhang, Zou]{Jianbing Lu, Yingnan Zhang, Hanlin Zou}

\address{Jianbing Lu, School of Mathematical Sciences, Zhejiang University, Hangzhou 310027, China}
\email{jianbinglu@zju.edu.cn}

\address{Yingnan Zhang, School of Statistics and Mathematics, Yunnan University of Finance and Economics, Kunming 650221, China}
\email{zyn2007@126.com}

\address{Hanlin Zou, School of Mathematics and Statistics, Yunnan University, Kunming 650091, China}
\email{zouhanlin@ynu.edu.cn}


\begin{document}

\begin{abstract}
A central problem in the study of generalized quadrangles is to classify finite generalized quadrangles satisfying certain symmetry conditions. It is known that an automorphism group of a finite thick generalized quadrangle $\mc{S}$ acting primitively on both the points and lines of $\mc{S}$ must be almost simple. In this paper, we initiate the study of finite generalized quadrangles admitting a point-primitive and line-primitive automorphism group with socle being a unitary group. We develop a group-theoretic tool to prove that the socle of such a group cannot be $\PSU(3,q)$ with $q\geq 3$.
\end{abstract}	

\maketitle

\section{Introduction}

A (finite) {\it generalized quadrangle} $(\mathrm{GQ})$ is a point-line incidence geometry $\mc{S}=(\mc{P},\mc{L},{\bf I})$ where $\mc{P}$ and $\mc{L}$ are disjoint sets of objects called points and lines, respectively, and {\bf I} is a symmetric point-line incidence relation satisfying the following axioms: (i) two distinct points are incident with at most one line, and (ii) given a line $l$ and a point $\a$ not incident with $l$, there is a unique point $\a'$ and line $l'$ such that $\a'$ is incident with $l$, and $l'$ is incident with both $\a$ and $\a'$. 
Generalized quadrangles can alternatively be defined as a point-line geometry whose incidence graph has diameter 4 and girth 8. In general, a finite {\it generalized $n$-gon (polygon)} is a point-line geometry whose incidence graph has diameter $n$ and girth $2n$ (see \cite[Lemma 1.3.6]{Van98} for example). We say that a generalized $n$-gon has {\it order} $(s,t)$ if there exist positive integers $s$ and $t$ such that each line is incident with $s+1$ points and each point is incident with $t+1$ lines. A generalized $n$-gon of order $(s,t)$ is called {\it thick} if $s>1$ and $t>1$. The Feit-Higman theorem \cite{FH1964} states that finite thick generalized $n$-gons exist if and only if $n=3,4,6$ or $8$.

Generalized $n$-gons were introduced by Tits \cite{Tits59} in an attempt to find geometric models for simple groups of Lie type. The polygons which arise from Lie type groups are called the {\it classical} generalized polygons. In particular, the classical triangles are the Desarguesian projective planes having $\mathrm{PSL}(3,q)$ as their automorphism groups, and the classical quadrangles are finite classical polar spaces admitting the groups $\mathrm{PSp}(4,q)$, $\mathrm{PSU}(4,q)$, or $\mathrm{PSU}(5,q)$. The classical hexagons are related to the $\mathrm{G}_{2}(q)$ and ${ }^{3} \mathrm{D}_{4}(q)$ groups, and the classical octagons correspond to the ${ }^{2} \mathrm{F}_{4}(2^{2k+1})$ groups. In all these cases, the groups act primitively on
both points and lines, and transitively on the {\it flags} (i.e., the incident point-line pairs) of the generalized polygon. The automorphism groups of the classical generalized polygons also act distance-transitively on the points of these geometries.

A natural question is: can any of these symmetry conditions characterize the classical generalized polygons? In \cite{BVM1994}, Buekenhout and Van Maldeghem proved that a finite generalized polygon with an automorphism group acting distance-transitively on the points is classical, dual classical, or the unique generalized quadrangle $\GQ(3,5)$ of order $(3,5)$.
For generalized quadrangles, Kantor \cite{Kantor} conjectured that a finite flag-transitive generalized quadrangle is classical, $\mathrm{GQ}(3,5)$ or the generalized quadrangle of order (15,17) arising from the Lunelli-Sce hyperoval up to duality. In the past few decades, a lot of work has been done towards this conjecture. For example, Bamberg, Li and Swartz classified antiflag-transitive generalized quadrangles in \cite{BLS2018} (an {\it antiflag} is a non-incident point-line pair) and locally 2-transitive generalized quadrangles in \cite{BLS2021}. This is the state-of-the-art regarding Kantor's conjecture. Another important problem is the following.
\begin{problem}\label{prob}
Classify finite thick generalized quadrangles admitting a point-primitive and line-primitive automorphism group.
\end{problem}
To this end, Bamberg et al. made a significant progress in \cite{BGMRS2012}.
\begin{theorem}[{\cite[Theorem 1.1]{BGMRS2012}}]\label{thm_plp}
An automorphism group $G$ of a finite thick generalized quadrangle $\mc{S}$ acting primitively on the points and lines is almost simple. In other words, the socle of $G$ is a non-abelian simple group.
\end{theorem}

In the same paper, they further eliminated almost simple groups whose socle is sporadic, and then those whose socle is alternating under an additional assumption that $G$ is flag-transitive. The major part left regarding this problem is therefore to consider the case when $G$ is almost simple of Lie type. In \cite{FL2022}, Feng and the first-named author considered the problem when $G$ is an almost simple linear group and proved that if the socle of $G$ is $\mathrm{PSL}(2,q)$ with $q\geq4$, then $\mathcal{S}$ is the symplectic quadrangle $W(2)$. In the present paper, we start investigating Problem \ref{prob} when the automorphism group is an almost simple unitary group and prove the following main result.

%
%

\begin{theorem}\label{main}
Let $G$ be an automorphism group of a finite thick generalized quadrangle $\mathcal{S}$.
If $G$ acts primitively on both points and lines of $\mc{S}$, then the socle of $G$ cannot be $\PSU(3,q)$ with $q\geq 3$.
\end{theorem}

The paper is organized as follows. In Section \ref{sec_prelim}, we present some preliminary results on generalized quadrangles and permutation groups. In Section \ref{sec_psuinfo}, we collect the maximal subgroups of almost simple groups with socle $\PSU(3,q)$. We also discuss the conjugacy classes and centralizers of elements of order 3 and 4 of $\PSU(3,q)$. With these preparations, the proof of Theorem \ref{main} is given in Section \ref{sec_pfmain}.

\section{Background on generalized quadrangles and permutation groups}\label{sec_prelim}

In this section we state some useful facts about both finite generalized quadrangles and permutation groups. Throughout this paper, we use $\mc{S}=(\mathcal{P}, \mathcal{L}, {\bf I})$ to denote a finite generalized quadrangle with point set $\mc{P}$, line set $\mc{L}$, and incidence relation ${\bf I}$.
 
\subsection{Parameters of generalized quadrangles}\label{GQ}

Recall that a (finite) generalized quadrangle $\mc{S}$ has order $(s,t)$ if every line contains $s+1$ points, and every point lies on $t+1$ lines. When $s=t$, we say $\mc{S}$ has order $s$. If we reverse the role of points and lines of $\mc{S}$, we obtain another generalized quadrangle called the {\it dual} of $\mc{S}$, and it has order $(t,s)$. The following two lemmas summarize a few basic results concerning the parameters of a generalized quadrangle.

\begin{lemma}[{\cite[1.2.1, 1.2.2, 1.2.3]{PT09}}]\label{para}
Let $\mathcal{S}$ be a generalized quadrangle of order $(s,t)$ with $v$ point and $b$ lines.
The following hold:
\begin{enumerate}[label={(\arabic*)}]
\item $v=(s+1)(st+1)$ and $b=(t+1)(st+1)$; \label{vb}
\item $s+t$ divides $st(s+1)(t+1)$; \label{divcond}
\item if $s>1$ and $t>1$, then $s\leq t^2$ and $t\leq s^2$. \label{st2}
\end{enumerate}
\end{lemma}

\begin{lemma}\label{newcond}
Let $\mathcal{S}$ be a thick generalized quadrangle of order $(s,t)$ with $v$ points and $b$ lines. The following hold:
\begin{enumerate}[label={(\arabic*)}]
\item $v>(s+1)^2(t+1)/2$, and $b>(t+1)^2(s+1)/2$; \label{cond1}
\item $v>({\frac{b}{v}})^4/2=(\frac{t+1}{s+1})^4/2$, and $b>(\frac{v}{b})^4/2=(\frac{s+1}{t+1})^4/2$. \label{cond2}
\end{enumerate}
\end{lemma}
\begin{proof}
\begin{enumerate}
\item Note that $(s-1)(t-1)>0$ as $\mc{S}$ is thick. Then $st+1>(s+1)(t+1)/2$ and the results follow from Lemma \ref{para}\ref{vb}.
\item By Lemma \ref{para}\ref{st2}, we have $(s+1)^2>t+1$. It then follows from Lemma \ref{newcond}\ref{cond1} that
$v>\frac{(s+1)^2(t+1)}{2}> (\frac{t+1}{s+1})^4/2$. Furthermore, Lemma \ref{para}\ref{vb} implies $\frac{v}{b}=\frac{s+1}{t+1}$. Thus $v>(\frac{b}{v})^4/2$. The other part is proved similarly.
\end{enumerate}
\end{proof}

\subsection{Permutation group theory}

We assume the reader is familiar with the basic notions in permutation group theory such as transitive actions, regular actions and primitive actions. Let $G$ be a group acting on a set $\Omega$. For $\a\in\Omega$ and $g\in G$, we use $\a^g$ to denote the image of $\a$ under the permutation corresponding to $g$, and use $G_\a$ to denote the stabilizer of $\a$ in $G$. The centralizer of $g$ in $G$ is denoted by $C_G(g)$ and the conjugacy class of $g$ in $G$ is denoted by $g^G$. More definitions and notation can be found in \cite{DM1996}.


Given a group $G$ acting on a set $\Omega$ and an element $g\in G$, define $\Fix(g)=\{x\in\Omega: x^g=x\}$. The next lemma gives a way of computing the number of fixed elements.
\begin{lemma}\cite[Lemma 2.5]{LS1991}\label{Fix}
Let $G$ be a finite transitive group acting on a set $\Omega$. Let $\a\in\Omega$ and $g\in G$. Then
\[|\Fix(g)|=\frac{|\Omega|\cdot|g^G\cap G_{\a}|}{|g^G|}.\]
\end{lemma}

This lemma enables us to prove the following.
\begin{lemma}\label{trans}
Let $G$ be a group acting transitively on a set $\Omega$ and let $g$ be a nontrivial element of $G$ fixing some $\a\in\Omega$. Then $C_G(g)$ is transitive on $\Fix(g)$ if and only if $g^G\cap G_{\a}=g^{G_{\a}}$.
\end{lemma}
\begin{proof}
Note that $C_G(g)$ stabilizes $\Fix(g)$ setwise and $C_{G_{\a}}(g)$ is the stabilizer of $\alpha$ in $C_G(g)$. Then $C_G(g)$ is transitive on $\Fix(g)$ if and only if $|\Fix(g)|=\frac{|C_G(g)|}{|C_{G_\alpha}(g)|}$. On the other hand, since $|g^G|=\frac{|G|}{|C_G(g)|}$, by Lemma \ref{Fix}, we have $|\Fix(g)|=\frac{|C_G(g)|\cdot |g^G\cap G_\alpha|}{|G_\alpha|}$. Thus $C_G(g)$ is transitive on $\Fix(g)$ if and only if $|g^G\cap G_\alpha|=\frac{|G_\alpha|}{|C_{G_\alpha}(g)|}=|g^{G_\alpha}|$, which is equivalent to $g^G\cap G_\alpha=g^{G_\alpha}$.
\end{proof}


\subsection{Substructures of a generalized quadrangle fixed by an automorphism}

Let $\mc{S}=(\mc{P},\mc{L},{\bf I})$ be a generalized quadrangle. We say that two points $x_1$ and $x_2$ are {\it collinear} if there is a line incident with both $x_1$ and $x_2$. Similarly, two lines $L_1$ and $L_2$ are {\it concurrent} if there is a point incident with both of them. We write $x_1\sim x_2$ if $x_1$ and $x_2$ are collinear or $x_1=x_2$, and dually, we write $L_1\sim L_2$ if $L_1$ and $L_2$ are concurrent or $L_1=L_2$.

A generalized quadrangle $\mc{S'}=(\mc{P}',\mc{L}',{\bf I}')$ is called a \emph{subquadrangle} of $\mc{S}=(\mc{P},\mc{L},{\bf I})$ if $\mc{P}'\subseteq \mc{P}$, $\mc{L}'\subseteq\mc{L}'$ and if ${\bf I}'$ is the restriction of ${\bf I}$ on $(\mc{P}'\times \mc{L}')\cup (\mc{L}'\times \mc{P}')$.

A \emph{grid} with parameters $(s_{1},s_{2})$ is a point-line incidence structure $(\mathcal{X}, \mathcal{B}, {\rm I})$ with
$$\mathcal{X}=\left\{x_{i, j}: 0 \leqslant i \leqslant s_{1}, 0 \leqslant j \leqslant s_{2}\right\}, \mathcal{B}=\left\{\ell_{0}, \ldots, \ell_{s_{1}}, \ell_{0}^{\prime}, \ldots, \ell_{s_{2}}^{\prime}\right\}$$
such that $x_{i, j}{\rm I} \ell_{k}$ if and only if $i=k$, and $x_{i, j} {\rm I} \ell_{k}^{\prime}$ if and only if $j=k$. A \emph{dual} \emph{grid} with parameters $(s_1,s_2)$ is the point-line dual of a grid with parameters $(s_2,s_1)$.

An {\it automorphism} (or a {\it collineation}) of $\mc{S}$ is a permutation of the points of $\mc{S}$ that preserves collinearity. Identifying each line with the points on it, an automorphism of $\mc{S}$ also permutes the lines of $\mc{S}$ and preserves concurrency. All automorphisms of $\mathcal{S}$ forms a group, which is called the full automorphism group and denoted by $\Aut(\mc{S})$. An automorphism group of $\mc{S}$ is a subgroup of $\Aut(\mc{S})$.

Let $g$ be an automorphism of $\mc{S}=(\mc{P},\mc{L},{\bf I})$. The substructure of $\mc{S}$ fixed by $g$ is an incidence structure $\mathcal{S}_{g}=\left(\mathcal{P}_{g}, \mathcal{L}_{g}, {\bf I}_{g}\right)$, where $\mathcal{P}_{g}$ and $\mathcal{L}_g$ are the set of fixed points and fixed lines of $g$ respectively, and ${\bf I}_g:={\bf I}\cap \big((\mathcal{P}_{g}\times \mathcal{L}_{g})\cup(\mathcal{L}_{g}\times \mathcal{P}_{g})\big)$. 

\begin{lemma}[{\cite[2.4.1]{PT09}}]\label{sub}
Let $g$ be an automorphism of a generalized quadrangle $\mathcal{S}=(\mathcal{P}, \mathcal{L}, {\bf I})$ of order $(s,t)$. The substructure $\mc{S}_g$ fixed by $g$ must be one of the following:
\begin{enumerate}
\item[(0)] $\mathcal{P}_{g}=\mathcal{L}_{g}=\varnothing$,
\item[(1)] $\mathcal{L}_{g}=\varnothing$, $\mathcal{P}_{g}$ is a nonempty set of pairwise noncollinear points,
\item[(1')] $\mathcal{P}_{g}=\varnothing$, $\mathcal{L}_{g}$ is a nonempty set of pairwise nonconcurrent lines,
\item[(2)] $\mathcal{P}_{g}$ contains a point $x$ such that $x \sim x'$ for every point $x' \in \mathcal{P}_{g}$, $\mathcal{L}_{g}$ is nonempty and each line of $\mathcal{L}_{g}$ is incident with $x$,
\item[(2')] $\mathcal{L}_{g}$ contains a line $L$ such that $L \sim L'$ for every line $L' \in \mathcal{L}_{g}$, $\mathcal{P}_{g}$ is nonempty and each point of $\mathcal{P}_{g}$ is incident with $L$,
\item[(3)] $\mathcal{S}_{g}$ is a grid with parameters $(s_1,s_2)$, $s_{1}<s_{2}$,
\item[(3')] $\mathcal{S}_{g}$ is a dual grid with parameters $(s_1,s_2)$, $s_{1}<s_{2}$,
\item[(4)] $\mathcal{S}_{g}$ is a subquadrangle of order $\left(s^{\prime}, t^{\prime}\right)$.
\end{enumerate}
\end{lemma}

\begin{lemma}\label{nonabelian}
Let $\cS$ be a thick generalized quadrangle of order $s$ admitting a regular automorphism group $G$ acting on the points. Then $G$ is nonabelian.
\end{lemma}
\begin{proof}
This follows immediately from \cite[Theorem 3.1]{Ghi1992}.
\end{proof}

Building on the above lemmas, we deduce the following corollary which will be used frequently in the proof of Theorem \ref{main}.
\begin{corollary}\label{centralizer}
Let $G$ be an automorphism group of a thick generalized quadrangle $\cS$ acting transitively on both points and lines. Suppose $G_\alpha\cap G_l\neq \varnothing$ for some point $\alpha$ and some line $l$. Let $g\in G_\alpha\cap G_l$ and $Y\leq C_G(g)$ be such that
\begin{enumerate}
\item $g^G\cap G_{\a}=g^{G_{\a}}$ and $g^G\cap G_l=g^{G_l}$;
\item $[Y:Y_\a]=[C_{G}(g):C_{G_\alpha}(g)]>4$ and $[Y:Y_l]=[C_{G}(g):C_{G_l}(g)]>4$.
\end{enumerate}
Then $Y$ is nonabelian.
\end{corollary}
\begin{proof}
Consider the substructure $\cS_g=(\cP_g,\cL_g,{\bf I}_g)$ fixed by $g$. By Lemma \ref{trans}, $C_G(g)$ is transitive on both $\mathcal{P}_g$ and $\mathcal{L}_g$. This implies that each point in $\cP_g$ is incident with a constant number of lines in $\cL_g$, and each line in $\cL_g$ is also incident with a constant number of points in $\cP_g$. Thus $\cS_g$ cannot have type (2), (2'), (3) and  (3') in Lemma \ref{sub}. We can further deduce from the transitivity of $C_G(g)$ that $|\cP_g|=[C_{G}(g):C_{G_\alpha}(g)]>4$ and $|\cL_g|=[C_{G}(g):C_{G_l}(g)]>4$. This means that $\cS_g$ cannot have type (0), (1) and (1'). Therefore $\cS_g$ is a subquadrangle. 

Since $[Y:Y_\a]=[C_{G}(g):C_{G_\alpha}(g)]$ and $[Y:Y_l]=[C_{G}(g):C_{G_l}(g)]$, then $Y$ is also transitive on both $\mc{P}_g$ and $\mc{L}_g$. Now suppose to the contrary that $Y$ is abelian. Then the kernel of $Y$ acting on $\cP_g$ is exactly $Y_\a$. So for any $h\in Y_\a$, it fixes every point in $\cP_g$ and thus fixes every line in $\cL_g$ since $\cS_g$ is a subquadrangle. So $Y_\a\subseteq Y_l$. Similarly, we have $Y_l\subseteq Y_\a$. Hence $Y_\a=Y_l$ which implies that $|\cP_g|=|\cL_g|>4$. Therefore $\cS_g$ is a subquadrangle of order $s'>1$ which admits a regular automorphism group $Y/Y_\a$. Finally, Lemma \ref{nonabelian} implies that $Y/Y_\a$ is nonabelian, which is contrary to the assumption that $Y$ is abelian. Therefore $Y$ is nonabelian.
\end{proof}

\section{Maximal subgroups and conjugacy classes of $\PSU(3,q)$}\label{sec_psuinfo}
Let $p$ be a prime, $f$ a positive integer and $q=p^f$. Let $\F_{q^2}$ be the finite field of order $q^2$. Using the automorphism $a\mapsto \bar{a}=a^q$ $(a\in\F_{q^2})$, we define
\[\SU(3,q):=\{A\in\SL(3,q^2)\mid AW\overline{A}^\top=W\},\]
where $\overline{A}=(\overline{a_{ij}})_{1\leq i,j\leq 3}$ and $W$ is any nonsingular matrix satisfying $W=\overline{W}^\top$. In fact, $W$ is the Gram matrix of a Hermitian form on the 3-dimensional vector space $\F_{q^2}^3$ over $\F_{q^2}$ with respect to a given basis. It is often convenient to choose $W$ to be the identity matrix $I$, or the anti-diagonal matrix $\antidiag(1,1,1)$. The groups $\SU(3,q)$ for different choices of $W$ are all conjugate in $\GL(3,q^2)$.

The projective special unitary group $\PSU(3,q)$ is defined to be $\SU(3,q)/Z(\SU(3,q))$, where $Z(\SU(3,q))$ is the center of $\SU(3,q)$. Let $~\hat{}~$ be the projection from $\SU(3,q)$ to $\PSU(3,q)$. We use $\hat{}H$ to represent the image of a subgroup $H$ of $\SU(3,q)$, and use $\hat{g}$ to represent the image of an element $g$ of $\SU(3,q)$.

\subsection{Maximal subgroups of almost simple groups with socle $\PSU(3,q)$}

Let $G$ be an almost simple group with socle $T= \PSU(3,q)$. In our study, we need a list of all maximal subgroups of $G$. Let $M$ be a maximal subgroup of $G$. The group $M_0:=T\cap M$ is usually maximal in $T$ and we call $M$ a {\it novelty} if $M_0$ is not maximal in $T$. 
By  \cite[Table 8.5 and Table 8.6]{lowdim}, the complete list of maximal subgroups of an almost simple group $G$ with socle $\PSU(3,q)$
is known, and in this case, there arise only four novelties.

\begin{lemma}\cite[Tables 8.5 and 8.6]{lowdim}\label{maximal}
Let $G$ be an almost simple group with socle $T=\PSU(3,q)$, $q=p^f\geq 3$, and let $M$ be a maximal subgroup of $G$ not containing $T$.
Then $M\cap T$ is isomorphic to one of the groups listed in Table \ref{tab_sub}. Moreover, $M\cap T$ is maximal in $T$ unless a novelty is indicated.
\end{lemma}

\begin{table}[h]\aboverulesep=0pt \belowrulesep=0pt
\setlength{\abovecaptionskip}{0cm}
\setlength{\belowcaptionskip}{0cm}
\caption{Maximal subgroups of $T=\PSU(3,q)$ and their indices in $T$}
\label{tab_sub}
\centering
\[\small
\arraycolsep=5pt
\begin{array}{llll}
\toprule
\text{}  & \text{Subgroup} & \text{Notes} & [T:M\cap T] \\
\midrule\midrule
&&&\\[-0.4cm]
\text{(a)} & \hat~[q]^{1+2}:(q^2-1) & & q^3+1\\[0.01cm]
\text{(b)} &\hat~\GU(2,q) && q^2(q^2-q+1)\\[0.01cm]
\text{(c)} &\hat~(q^2-q+1):3 & q\neq 3,5  & q^3(q^2-1)(q+1)/3\\
\text{(c')}&&\text{novelty if }q=5 \\
\text{(d)}&\hat~(q+1)^2:\sym_3  & q\neq 5& q^3(q-1)(q^2-q+1)/6\\
\text{(d')}&& \text{novelty if }q=5\\
\text{(e)}&\SO(3,q)&  q\geq7, q \text{ odd} & q^2(q^3+1)/\gcd(3,q+1)\\[0.1cm]
\text{(f)}&\PSU(3,q_0). \gcd(3,\frac{q+1}{q_0+1}) & q=q_0^r, r \text{ odd and prime} & \frac{1}{c}\!\cdot\!\frac{q^3(q^2-1)(q^3+1)}{q_0^3(q_0^2-1)(q_0^3+1)},  c\!=\!\gcd(3,\frac{q+1}{q_0+1})\\[0.2cm]
\text{(g)}&\hat~3_{+}^{1+2}:\rQ_8.\frac{\gcd(9,q+1)}{3}&  p\!=\!q\!\equiv\!2\hspace{-0.2cm}\pmod 3, q\geq 11 & \frac{q^3(q^2-1)(q^3+1)}{72\cdot\gcd(9,q+1)}\\
\text{(g')}&&  \text{novelty if }q=5\\
\text{(h)} &\PSL(2,7) & q\!\neq \!5, p\!=\!q\!\equiv\!3,\!5,\!6\hspace{-0.2cm}\pmod 7 &  \frac{q^3(q^2-1)(q^3+1)}{168\cdot\gcd(3,q+1)}\\
\text{(h')}&& \text{novelty if }q=5\\
\text{(i)}&\A_6 & p\!=\!q\!\equiv \!11,14\hspace{-0.2cm}\pmod {15}& \frac{q^3(q^2-1)(q^3+1)}{360\cdot\gcd(3,q+1)}\\
\text{(j)}&\A_6. 2_3 & q=5\\
\text{(k)}&\A_7 & q=5\\
\bottomrule
\end{array}
\]
\end{table}

\subsection{Conjugacy classes and centralizers of elements in $\PSU(3,q)$}

We shall make use of the conjugacy classes and the centralizers of elements of order 3 and 4 of $\PSU(3,q)$ in the proof of the main theorem. The conjugacy classes and the orders of the centralizers have been worked out in \cite[Table 2]{SF1973}. The structures of the centralizers are determined in the next two lemmas. All these results are summarized in Table \ref{tab_ord3} and Table \ref{tab_ord4} for future reference. 
Note that the representatives of different conjugacy classes may belong to distinct $\PSU(3,q)$'s defined by different Gram matrices. We will indicate which Gram matrix is used in the proofs when needed.

In the sequel, we fix the following notation. Let $\F_{q^2}^*=\<\w\>$ be the multiplicative group of the field $\F_{q^2}$. Denote by $\zeta_n$ an element of order $n$ in $\F_{q^2}^*$. 

\begin{table}[h]\aboverulesep=0pt \belowrulesep=0pt
\setlength{\abovecaptionskip}{0.2cm}
\setlength{\belowcaptionskip}{0cm}
\caption{{\spaceskip=0.27em\relax  Conjugacy classes and centralizers of elements of order 3 of $\PSU(3,q)$}}
$T=\PSU(3,q)$, $\rho=\w^{q-1}$, $x_0+x_0^q=0$
\centering
\[\footnotesize
\begin{array}{>{\hfil$}p{2.5cm}<{$\hfil}|c|>{\hfil$}p{4.4cm}<{$\hfil}>{\hfil$}p{2.6cm}<{$\hfil}>{\hfil$}p{3cm}<{$\hfil}}
\toprule
& & &\\[-0.4cm]
q & \text{\makecell[c]{\# of\\ classes}} & \text{Representative~} g&  \text{Gram matrix~}W & |C_T(\widehat{g})|\\[0.1cm]
\midrule\midrule
& & &\\[-0.3cm]
\multirow{6}{*}{$q\equiv 0\pmod{3}$}   &  \multirow{6}{*}{$2$} & {\footnotesize\begin{pmatrix}1& 0&0\\1&1&0\\ 0&1 &1\end{pmatrix}} & {\footnotesize\begin{pmatrix}0& 0&-1\\0&1&1\\ -1&1 &0\end{pmatrix}}& q^2\\[0.5cm]
\cmidrule{3-5}
& & &\\[-0.3cm]
& &{\footnotesize\begin{pmatrix}1& 0&0\\1&1&0\\ 0&0 &1\end{pmatrix}} & {\footnotesize\begin{pmatrix}0& x_0&0\\ \overline{x_0}&0&0\\ 0&0 &1\end{pmatrix}} & q^3(q+1) \\[0.5cm]
\midrule
& & &\\[-0.3cm]
q\equiv 1\pmod{3} & 1 &{\footnotesize\begin{pmatrix}1&0&0\\ 0 &\zeta_3 &0 \\ 0 &0 &\zeta_3^2\end{pmatrix}} & {\footnotesize\begin{pmatrix}1&0 &0\\ 0&1&0\\ 0&0 &1\end{pmatrix}}& q^2-1\\[0.5cm]
\midrule
& & &\\[-0.3cm]
\makecell[c]{q\equiv 2\pmod{3}\\ \text{and} \\9\nmid q+1}& 1&{\footnotesize\begin{pmatrix}1&0&0\\ 0 &\zeta_3 &0 \\ 0 &0 &\zeta_3^2\end{pmatrix}} & {\footnotesize\begin{pmatrix}1&0 &0\\ 0&1&0\\ 0&0 &1\end{pmatrix}}&(q+1)^2\\[0.5cm]
\midrule
& & &\\[-0.3cm]
\multirow{7}{*}{\makecell[c]{$q\equiv 2 \pmod{3}$\\ \text{and} \\ $9\mid q+1$}}  & \multirow{7}{*}{$3$}&  {\footnotesize\begin{pmatrix}1&0&0\\ 0 &\zeta_3 &0 \\ 0 &0 &\zeta_3^2\end{pmatrix}} & {\footnotesize\begin{pmatrix}1&0 &0\\ 0&1&0\\ 0&0 &1\end{pmatrix}}&(q+1)^2\\[0.6cm]
\cmidrule{3-5}
& & &\\[-0.3cm]
& & {\footnotesize\begin{pmatrix}\rho^{\frac{k(q+1)}{9}}&\hspace{-0.2cm}0&\hspace{-0.2cm}0\\ 
0 &\hspace{-0.2cm}\rho^{\frac{k(q+1)}{9}} &\hspace{-0.2cm}0 \\ 
0 &\hspace{-0.2cm}0 &\hspace{-0.2cm}\rho^{-\frac{2k(q+1)}{9}}\end{pmatrix} }& {\footnotesize\begin{pmatrix}1&0 &0\\ 0&1&0\\ 0&0 &1\end{pmatrix}}&\displaystyle\frac{q(q+1)^2(q-1)}{3}\\[0.8cm]
& &{\text{\footnotesize $k\in\{1,2\}$}}&\\[0.1cm]
\bottomrule
\end{array}
\]
\label{tab_ord3}
\end{table}

\begin{lemma}\label{order3}
Suppose $q\equiv 2\pmod{3}$ and $q>2$. Let $X=\{A\in\SL(3,q^2)\mid A\overline{A}^\top=I\}$ and 
let $T=\hat{} X=\PSU(3,q)$. Assume $g=\diag(1,\zeta_3,\zeta_3^2)$. Then $C_T(\widehat{g})=(\C_{(q+1)/3}\times \C_{q+1}):\C_3$, which is nonabelian.
\end{lemma}
\begin{proof} 
Define \[c_1={\scriptsize\begin{pmatrix}\w^{q-1}&0&0\\0&1&0\\0&0&\w^{1-q}\end{pmatrix}}, c_2={\scriptsize\begin{pmatrix}1&0&0\\ 0&\w^{q-1}&0\\ 0&0&\w^{1-q}\end{pmatrix}},\text{ and }c_3={\scriptsize\begin{pmatrix}0&1&0\\0&0&1\\1&0&0\end{pmatrix}}.\]
It is easy to check that $g^{c_1}=g, g^{c_2}=g$, and $g^{c_3}=\zeta_3^2g$. Furthermore, we have $c_1c_2=c_2c_1$, $\<c_1\>\cap\<c_2\>=\<I\>$, and $c_1^{c_3}, c_2^{c_3}\in\<c_1,c_2\>=\<c_1c_2\>\times\<c_2\>$. Therefore \[C_T(\widehat{g})=(\<\widehat{c_1c_2}\>\times \<\widehat{c_2}\>): \<\widehat{c_3}\>\cong (\C_{(q+1)/3}\times \C_{q+1}):\C_3\] since $|C_T(\widehat{g})|=(q+1)^2$ according to Table \ref{tab_ord3}. Finally, a direct computation shows that $\widehat{c_2}\widehat{c_3}\neq \widehat{c_3}\widehat{c_2}$. Therefore $C_T(\widehat{g})$ is nonabelian.
\end{proof}

\begin{table}[h]\aboverulesep=0pt \belowrulesep=0pt
\setlength{\abovecaptionskip}{0.2cm}
\setlength{\belowcaptionskip}{0cm}
\caption{{\spaceskip=0.27em\relax Conjugacy classes and centralizers of elements of order 4 of $\PSU(3,q)$}}
\label{tab_ord4}
$T=\PSU(3,q)$, $d=\gcd(3,q+1)$, $\rho=\w^{q-1}$
\centering
\[\footnotesize
\arraycolsep=5pt
\begin{array}{c|c|cccc}
\toprule
& & &\\[-0.4cm]
q & \text{\makecell[c]{\# of\\ classes}} & {\rm Representative~} g & {\rm Gram~matrix~}W &|C_T(\widehat{g})|&C_T(\widehat{g})\\[0.1cm]
\midrule\midrule
& & &\\[-0.3cm]
q\equiv 1\!\!\!\! \pmod{4}\!\!\!\! &  1 & {\footnotesize\begin{pmatrix}1& 0&0\\0&\zeta_4&0\\ 0&0 &\zeta_4^3\end{pmatrix}} & {\footnotesize\begin{pmatrix}1& 0&0\\0&0&1\\ 0&1 &0\end{pmatrix}} & \displaystyle\frac{q^2-1}{d} & \C_{(q^2-1)/d}\\[0.5cm]
\midrule
& & &\\[-0.3cm]
\multirow{9}{*}{$q\equiv 3 \!\!\!\!\pmod{4}$\!\!\!\! } & \multirow{9}{*}{$3$} &{\footnotesize\begin{pmatrix}\rho^k&0&0\\ 0 &\rho^k &0 \\ 0 &0 &\rho^{-2k}\end{pmatrix}}& {\footnotesize\begin{pmatrix}1&0 &0\\ 0&1&0\\ 0&0 &1\end{pmatrix}}&\displaystyle\frac{q(q\!+\!1)^2(q\!-\!1)}{d} & \hat{~}\GU(2,q)\\[0.5cm]
& & \text{\footnotesize $k\in\left\{\frac{q+1}{4d},\frac{3(q+1)}{4d}\right\}$}& \\[0.3cm]
\cmidrule{3-6}
& & &\\[-0.2cm]
& &{\footnotesize\begin{pmatrix}\rho^k&0&0\\ 0 &\rho^\ell &0 \\ 0 &0 &\rho^m\end{pmatrix}} &{\footnotesize\begin{pmatrix}1&0 &0\\ 0&1&0\\ 0&0 &1\end{pmatrix}}& \displaystyle\frac{(q+1)^2}{d} &\C_{(q+1)/d}\!\times\! \C_{q+1}\\[0.6cm]
& & \multicolumn{3}{l}{\text{\footnotesize$(k,\ell,m)\!=\!\begin{cases}(\frac{q+1}{12},\frac{q+1}{3},\frac{7(q+1)}{12}),&\!\!\!\!\!\text{if } q\equiv 2\!\!\!\!\pmod{3}\\
(\frac{q+1}{4},\frac{3(q+1)}{4}, q\!+\!1),&\!\!\!\!\!\text{if } q\not\equiv 2\!\!\!\!\pmod{3}\end{cases}$}} & \\[0.4cm]
\midrule
& & &\\[-0.3cm]
q=2^f, f\geq 1 & d & {\footnotesize \begin{pmatrix}1&0&0\\ z_\ell &1 &0 \\ 0 & \overline{z_\ell} &1\end{pmatrix}} &{\footnotesize\begin{pmatrix}0&0 &1\\ 0&1&1\\ 1&1 &0\end{pmatrix}}& q^2 &\C_4^f\\[0.6cm]
&&\ell\in\{1,\ldots,d\},\\
&& \multicolumn{2}{l}{z_\ell\in\{\rho^k+1\mid k\equiv \ell\!\!\!\!\pmod{3}\}}\\[0.1cm]
\bottomrule
\end{array}
\]
\flushleft
\vspace{1em}

{\it Remarks.} 

(i) Everything in Table \ref{tab_ord4} was taken from \cite[Table 2]{SF1973} except the last column whose entries will be determined in Lemma \ref{order4}.

(ii) When $q$ is even, the representatives given in \cite[Table 2]{SF1973} are actually: ${\scriptsize \begin{pmatrix}1&0&0\\ \theta^{\ell} &1 &0 \\ 0 &\theta^{\ell} &1\end{pmatrix}}$ where $\theta\in\F_{q^2}\setminus\F_q$, $\theta^3\neq 1$ and $\ell\in\{0,\ldots,d-1\}$. For our purposes, we use the elements described in Table \ref{tab_ord4} as representatives. It can be verified by direct computation that these elements are pairwise non-conjugate when $d=3$.
\end{table}

\begin{lemma}\label{order4}
Let $g$ and $W$ be as specified in Table \ref{tab_ord4}. Let $X=\{A\in\SL(3,q^2)\mid AW\overline{A}^\top=W\}$ and let $T=\hat~X=\PSU(3,q)$. Then the following hold.
\begin{enumerate}[label={(\arabic*)}]
\item If $q\equiv 1\pmod{4}$, then $C_T(\widehat{g})=\C_{(q^2-1)/d}$.\label{ord4.1}
\item If $q\equiv 3\pmod{4}$ and $g=\diag(\rho^k,\rho^k,\rho^{-2k})$, then $C_T(\widehat{g})=\hat{~}\GU(2,q)$.
\item If $q\equiv 3\pmod{4}$ and $g=\diag(\rho^k,\rho^\ell,\rho^m)$, then $C_T(\widehat{g})=\C_{(q+1)/d}\times \C_{q+1}$.
\item If $q=2^f$ with $f\geq 1$, then $C_T(\widehat{g})=\C_4^f$.\label{ord4.4}
\end{enumerate}
\end{lemma}
\begin{proof}The orders of the centralizers were determined in \cite[Table 2]{SF1973} and we have included them in Table \ref{tab_ord4}. In the following, we determine their structures.
\begin{enumerate}
\item Let $c={\scriptsize\begin{pmatrix}\w^{q-1} &0&0\\0&\w^{-q}&0\\0&0&\w\end{pmatrix}}$. Note that $g^{c}=g$ and the order of $\widehat{c}$ is $(q^2-1)/d$ which is equal to $|C_T(\widehat{g})|$. Therefore $C_T(\widehat{g})=\<\widehat{c}\,\>=\C_{(q^2-1)/d}$.
\item By Table \ref{tab_sub}, we see that there is a unique maximal subgroup of $T$ whose order is a multiple of $|C_T(\widehat{g})|$, namely $\hat~\GU(2,q)$. Moreover, since $|\hat~\GU(2,q)|=|C_T(\widehat{g})|$, we conclude that $C_T(\widehat{g})=\hat~\GU(2,q)$.

\item 
Take $c_1$ and $c_2$ to be the matrices in the proof of Lemma \ref{order3}. Then it can be proved similarly that $C_T(\hat{g})=\<\widehat{c_1c_2}\>\times \<\widehat{c_2}\>\cong \C_{(q+1)/d}\times \C_{q+1}$.

\item Take $g={\scriptsize \begin{pmatrix}1&0&0\\ z_\ell &1 &0 \\ 0 &\overline{z_\ell} &1\end{pmatrix}}$ where $\ell\in\{1,\ldots,d\}$ and $z_\ell\in\{\rho^k+1\mid k\equiv \ell\pmod{3}\}$. Note that $g\in X=\{A\in \SL(3,q^2) \mid AW\overline{A}^\top=W\}$ with $W={\scriptsize\begin{pmatrix}0&0 &1\\ 0&1&1\\ 1&1&0\end{pmatrix}}$. 
Define
\[H=\left\{\left.h_{a,b}:={\footnotesize\begin{pmatrix} 1&0&0\\ a &1&0\\ b& a & 1\end{pmatrix} }\right\vert \makecell[l]{a,b\in\F_{q^2}, \overline{z_\ell}a+z_\ell \overline{a}=0, \\[0.2cm]  b+\overline{b}+a+\overline{a}+a\overline{a}=0}\right\}.\]
It is readily seen that $H$ is an abelian subgroup of $X$ of order $q^2$. Moreover, each $h_{0,b}\in H$ with $b\neq 0$ has order $2$, and each $h_{a,b}\in H$ with $a\neq 0$ has order 4. This shows that the number of involutions of $H$ is exactly $q-1$. On the other hand, the fundamental theorem of finitely generated abelian groups implies $H=\C_{2}^i\times \C_{4}^j$ with $i+2j=2f$. It follows that the number of involutions of $H$ is $2^{i+j}-1\geq 2^f-1=q-1$. Thus $i+j=f$. So $j=f$ and $i=0$, which imply $H=\C_4^f$. Finally, noticing that $H\leq C_X(g)$, $H\cap Z(X)=\<I\>$ and $|C_T(\widehat{g})|=q^2$, we conclude that $C_T(\widehat{g})\cong H=\C_4^f$.
\end{enumerate}
\end{proof}

\section{Proof of Theorem \ref{main}}\label{sec_pfmain}
Our goal in this section is to prove Theorem \ref{main}.
Let $\mathcal{S}=(\mathcal{P}, \mathcal{L}, {\bf I} )$ be a finite thick generalized quadrangle admitting an automorphism group $G$ acting primitively on both points and lines. Let $T$ be the socle of $G$ with $T=\PSU(3,q)$, $q=p^f\geq 3$ and $p$ prime. For $\a\in \mathcal{P}$ and $l\in \mathcal{L}$, we have 
\[T_{\a}=G_{\a}\cap T,\text{ and }T_l=G_l\cap T.\]
Since $G$ is primitive on $\mathcal{P}$ and $\mathcal{L}$, and $T$ is normal in $G$, then $T$ is transitive on $\mathcal{P}$ and $\mathcal{L}$.
Let $v$ and $b$ be the number of points and the number of lines of $\mc{S}$, respectively. We have
\begin{align}
v=&[G:G_{\a}]=[T:T_{\a}],\label{point}\\
b=&[G:G_l]=[T:T_l],\label{line}
\end{align}
where $|T|=q^3(q^2-1)(q^3+1)/\gcd(3,q+1)$.
Since $G_\alpha$ (and $G_l$) is maximal in $G$, then $T_\alpha$ (and $T_l$) has to be isomorphic to one of the groups listed in Table \ref{tab_sub}. We will show in the following lemmas that, however, $T_\alpha$ cannot be isomorphic to any of them and therefore complete the proof of the main theorem.

Before going into details, let us briefly explain how we exclude all the possibilities of $T_\alpha$. The proof will be done by contradiction. By \eqref{point} and \eqref{line}, one can read off $v$ and $b$ directly from the last column of Table \ref{tab_sub}. In many cases, we are able to deduce an arithmetic contradiction by solving the system of equations $(1+s)(1+st)=v$ and $(1+t)(1+st)=b$. 
When it is harder to find a contradiction arithmetically, we consider the fixed substructure of certain element of $T$ and use group theoretic properties such as Corollary \ref{centralizer}. In this way, the remaining possibilities of $T_\alpha$ are ruled out.

For the sake of convenience, we shall say $T_\alpha$ (or $T_l$) has {\it type} (x) if $T_\alpha$ (or $T_l$) is isomorphic to the subgroup in row (x) of Table \ref{tab_sub}. For example, $T_\alpha$ has type (a) if $T_\alpha\cong \hat~[q]^{1+2}:(q^2-1)$.

First of all, it can be easily shown that $T_\a$ does not have certain types.
 
\begin{lemma}\label{a}
The subgroup $T_\a$ cannot have type (a), i.e., $T_\a\not\cong \hat~[q]^{1+2}:(q^2-1)$.
\end{lemma}
\begin{proof}
Suppose to the contrary that $T_\a\cong\hat~[q]^{1+2}:(q^2-1)$. By \cite[Table 2.3]{lowdim}, we see that $T_\a$ is the stabilizer of a 1-dimensional totally isotropic subspace of the natural module $V=\F_{q^2}^3$. Thus we may identify the point set of $\mc{S}$ with the set $V_1$ consisting of all 1-dimensional totally isotropic subspaces of $V$. It is well-known that $T=\PSU(3,q)$ is 2-transitive on $V_1$; see for example \cite[Section 7.7]{DM1996}. However, as an automorphism group of a generalized quadrangle, $T$ cannot be 2-transitive on the points because it cannot map two collinear points to two noncollinear points. Therefore $T_\a\not\cong\hat~[q]^{1+2}:(q^2-1)$.
\end{proof}

\begin{lemma}\label{q=5}
The subgroup $T_\alpha$ cannot have type (c'), (d'), (g'), (h'), (j) or (k).
\end{lemma}
\begin{proof}
Since $\mathcal{S}$ is thick, by Lemma \ref{para}\ref{vb}, we have $v\geq 15$ and $b\geq 15$. If $T_\alpha$ has type (c'),(d'),(g'),(h'),(j) or (k), then $q=5$. This implies $G$ is one of $\PSU(3,5)$, $\PGU(3,5)$, $\PSU(3,5).2$ and $\PGU(3,5).2$. In all cases, the indices of the maximal subgroups of $G$ can be computed by Mamga \cite{Magma}, so we can use $v=[G:G_\a]$ and $b=[G:G_l]$ to get $v,b\in$ $\{6000$, $1750$, $750$, $525$, $175$, $126$, $50\}$. However, for these values of $v$ and $b$, the system of equations $v=(s+1)(st+1)$ and $b=(t+1)(st+1)$ has no integer solutions for $s$ and $t$. Therefore $q$ cannot be 5 and the statement is proved.
\end{proof}

%

To eliminate the remaining possibilities of $T_\a$, we shall proceed in two cases according as $T_\a$ and $T_l$ have different types or not. The following is an immediate corollary of Lemma \ref{newcond}\ref{cond2} and it will exclude many combinations of types of $T_\a$ and $T_l$ very quickly. 

\begin{lemma}\label{lem_boundM1}
With the above notation, we have $\frac{1}{2^{1/5}}[T:T_\a]^{4/5}<[T:T_l]<2^{1/4}[T:T_\a]^{5/4}$.
\end{lemma}
\begin{proof}
The statement follows directly from Lemma \ref{newcond}\ref{cond2} by using $v=[T:T_\a]$ and $b=[T:T_l]$.
\end{proof}

Let $n$ be an integer and $\mu$ a prime. We will use $(n)_\mu$ to denote the {\it $\mu$-part} of $n$ which is the highest power of $\mu$ dividing $n$.

\subsection{The case when $T_\a$ and $T_l$ have different types}

By duality, we assume that the type of $T_l$ is always behind that of $T_\a$ in alphabetical order.

\begin{lemma}\label{b}
The subgroup $T_{\a}$ cannot have type (b), i.e., $T_\a\not\cong \hat~\GU(2,q)$.
\end{lemma}
\begin{proof}
Suppose to the contrary that $T_{\a}\cong\hat~\GU(2,q)$. First of all, one can obtain $[T:T_\a]$ and $[T:T_l]$ from the last column of Table \ref{tab_sub} as expressions of $q$ with certain restrictions. Then by the bounds on $[T:T_l]$ given in Lemma \ref{lem_boundM1}, we deduce that $T_l$ cannot have type (c), (g) or (i). Similarly, if $T_l\cong\hat~(q+1)^2:\sym_3$, then $q\leq8$; if $T_l\cong\PSU(3,q_0). \gcd(3,\frac{q+1}{q_0+1})$, then $r=3, q_{0}=2$; and if $T_l\cong\PSL(2,7)$, then $q=3$. So $(v,b)$ can take only a finite number of values, and for these values of $v$ and $b$, it is easy to check that there no feasible $(s,t)$ pairs that satisfy both $v=(s+1)(st+1)$ and $b=(t+1)(st+1)$. Thus $T_l$ cannot have type (d), (f) or (h).

Now we consider the remaining possibility, i.e., $T_l\cong\SO(3,q)$ with $q\geq7$ odd. When $q\leq 15$, it is easy to verify that there are no feasible $(s,t)$ pairs satisfying $v=(s+1)(st+1)$ and $b=(t+1)(st+1)$. In the rest of the proof, we suppose that $q>15$. By \eqref{point}, \eqref{line} and Table \ref{tab_sub}, we have \[v=[T:T_\a]=q^2(q^2-q+1)\text{ and }b=[T:T_l]=q^2(q^3+1)/d,\] where $d=\gcd(3,q+1)$. Thus $\frac{t+1}{s+1}=\frac{v}{b}=\frac{q+1}{d}$. If $d=1$, then $t+1=(q+1)(s+1)$. By Lemma \ref{para}\ref{st2}, we have $(q+1)(s+1)=t+1<(s+1)^2$ which implies $s>q$ and $t>q^2+2q$. It follows that $(s+1)(st+1)>q^4>v$, a contradiction. Hence $d=3$ and 
\begin{equation}\label{4.4t}
t=(q+1)(s+1)/3-1.
\end{equation}

If $s+1\geq 2q$, then
\begin{align*}
v=(s+1)(st+1)\geq &2q\left((2q-1)\Big(\frac{q+1}{3}\cdot 2q-1\Big)+1\right)>q^{4}>q^2(q^2-q+1)=v,
\end{align*} a contradiction. Hence $s+1<2q$. 

Similarly, we can deduce that $s+1>q$ and  write 
\begin{equation}\label{4.4s}
s+1=q+a_0,
\end{equation}where $0\leq a_0<q$.

Since $q$ is odd, then $q^2(q^2-q+1)=v=(s+1)(st+1)$ implies that $s+1$ is odd. So $\gcd(s+1,st+1)=\gcd(s+1,2)=1$. Moreover, since $q<s+1<2q$ and $q^2$ divides $v$ which equals $(s+1)(st+1)$, we have $q^2\mid st+1$. By \eqref{4.4t} and \eqref{4.4s}, we have
\begin{equation}\label{4.4st1}
0\equiv 3(st+1)\equiv6-4a_0+a_0^2+(-4+a_0+a_0^2)q  \pmod {q^{2}}.
\end{equation} 
We assume 
\begin{equation}\label{b1}
6-4a_0+a_0^2=mq
\end{equation}
for some integer $m$. Since $0\leq a_0<q$, we have $0<m<q$. Using \eqref{b1}, the congruence equation \eqref{4.4st1} can be reduced to $0\equiv (m-4+a_0+a_0^2)q\equiv (m+5a_0-10+mq)q\pmod{q^2}$. Thus $q\mid m+5a_0-10$. Set
\begin{equation}\label{b2}
m+5a_0-10=jq,
\end{equation}
for some integer $j$. We have $1\leq j\leq5$ since $0\leq a_0,m<q$.

By \eqref{4.4t} and \eqref{4.4s}, we have 
\begin{align*}
3q^2(q^2-q+1)=&3v=3(s+1)(st+1)\\
=&(a_0 + q) (6 - 4 a_0 + a_0^2 - 4 q + a_0 q + a_0^2 q + 2 a_0 q^2 + q^3).
\end{align*}
By \eqref{b1} and \eqref{b2}, we have $a_0^2=(4-5q)a_0+jq^2+10q-6$. Then the above equation can be further reduced to 
\[3q^2(q^2-q+1)=(18-2a_0+a_0j)q^2+(-20+13a_0+j+a_0j)q^3+(1-2j)q^4,\]
which implies $3q^2\equiv (18-2a_0+a_0j)q^2\pmod{q^3}$. Thus \[q\mid 15-2a_0+a_0j.\] Since $q>15$, $1\leq j\leq 5$ and $0\leq a_0<q$, we deduce that $(j, a_{0}, s+1)=(1, 15, q+15), (3, q-15, 2q-15), \left(4, \frac{q-15}{2}, \frac{3q-15}{2}\right)$, or $(5, q-5, 2q-5)$. We check by computer that there are no feasible $(s,t)$ pairs that satisfy $v=(s+1)(st+1)$ and $b=(t+1)(st+1)$ for these cases. Thus the last possibility of $T_l$ is excluded and this completes the proof.
\end{proof}

\begin{lemma}\label{c}
The subgroup $T_{\a}$ cannot have type (c), i.e., $T_\a\not\cong\hat~(q^2-q+1)\!:\!3$ with $q\neq 3,5$.
\end{lemma}
\begin{proof}
Suppose to the contrary that $T_{\a}\cong\hat~(q^2-q+1):3$ with $q\neq 3,5$. By \eqref{point} and Table \ref{tab_sub}, we see that \[v=[T:T_\a]=q^{3}(q^{2}-1)(q+1)/3.\] 
Using the bounds on $[T:T_l]$ in Lemma \ref{lem_boundM1}, we deduce that $q\leq 38096$ if $T_l\cong\hat~3_{+}^{1+2}:\rQ_8.\frac{\gcd(9,q+1)}{3}$; $q\leq 31$ if $T_l\cong\PSL(2,7)$; and $q\leq 49$ if $T_l\cong\A_6$.  We verify by computer that there are no feasible $(s,t)$ pairs that satisfy $v=(s+1)(st+1)$ and $b=(t+1)(st+1)$ in these cases. 

The possibilities that $T_l\cong\hat~(q+1)^2:\sym_3$ or $T_l\cong\SO(3,q)$ with $q\geq7$ odd can be excluded in a similar fashion, so we only give the details for exlucding $T_l\cong\hat~(q+1)^2:\sym_3$ here. In this case, by \eqref{line} and Table \ref{tab_sub}, we have $b=[T:T_l]=q^3(q-1)(q^2-q+1)/6$, so $\frac{s+1}{t+1}=\frac{v}{b}=\frac{2(q+1)^2}{q^2-q+1}$. Note that $\gcd\left(2(q+1)^2, q^2-q+1\right)=\gcd(3,q+1)=:d$. Then there exists a positive integer $k$ such that 
\[s+1=2(q+1)^2k/d\text{ and }t+1=(q^2-q+1)k/d.\] If $d=3$ and $k=1$, then $v=(s+1)(st+1)$ is reduced to $(q+1)^3 (5 q^3- 18 q^2+ 36 q-22)=0$ which has no positive integer solutions for $q$. If $d\neq 3$ or $k>1$, then by Lemma \ref{newcond}\ref{cond1}, we have $q^3(q^2-1)(q+1)/3=v>(s+1)^2(t+1)/2\geq 16(q+1)^3(q^3+1)/27$. But this is incorrect for any $q$. Therefore $T_l\not\cong\hat~(q+1)^2:\sym_3$.

Now we consider the remaining case, i.e., $T_l\cong\PSU(3,q_0). \gcd(3,\frac{q+1}{q_0+1})$, where $q=q_0^r$, $r$ odd and prime. We see from Table \ref{tab_sub} that 
\begin{equation}\label{4.5b}
b=\frac{q^3(q^2-1)(q^3+1)}{c\cdot q_0^3(q_0^2-1)(q_0^3+1)},
\end{equation} 
where $c=\gcd(3,\frac{q+1}{q_0+1})$. 
By Lemma \ref{newcond}\ref{cond2}, we have
\begin{align*}
q^{3}(q^{2}-1)(q+1)/3=v>&\frac{1}{2}\left(\frac{b}{v}\right)^4=\frac{1}{2}\left(\frac{3}{c}\cdot\frac{q^{2}-q+1}{q_0^3(q_0^2-1)(q_0^3+1)}\right)^4\\
\geq &\frac{1}{2}\left(\frac{3}{c}\cdot\frac{q^{2}-q+1}{q_0^3\cdot q_0^5}\right)^4\geq \frac{1}{2}\left(\frac{q^{2}-q+1}{q_0^8}\right)^4,
\end{align*}
which implies that $q^{3}(q^{2}-1)(q+1)/3\cdot q_{0}^{32}>(q^2-q+1)^4/2$. Since $q(q^2-1)(q+1)/3<(q^2-q+1)^2/2$, we deduce that $q_0^{2r+32}> (q_0^{2r}-q_0^r+1)^2$. If $r\geq 17$, then we have $q_0^{4r}\geq 4q_0^{2r+32}>4 (q_0^{2r}-q_0^r+1)^2$, which is impossible. Therefore, $r=3,5,7,11,13$.

Suppose that $r=3$. Then $c=\gcd(3,q_0^2-q_0+1)=\gcd(3,q_0+1)$ and 
\begin{equation}\label{cfr3st}
\frac{t+1}{s+1}=\frac{b}{v}=\frac{q_{0}^6-q_{0}^{3}+1}{(c/3)\cdot q_0^3(q_0^2-1)(q_0^3+1)}.
\end{equation}
Note that $\gcd\left(q_0^6-q_0^3+1,q_0^2-1\right)=\gcd(q_0^6-q_0^3+1,q_0^3+1)=\gcd(3,q_0+1)=c$, and $q_0^6-q_0^3+1\equiv 1$ or $3\pmod{9}$. Hence $\gcd\left(\frac{q_0^6-q_0^3+1}{c}, \frac{q_0^3(q_0^2-1)(q_0^3+1)}{3}\right)=1$. By \eqref{cfr3st}, we can assume
\begin{align*}
t+1=\frac{q_0^6-q_0^3+1}{c}k, \text{~and~}s+1=\frac{q_0^3(q_0^2-1)(q_0^3+1)}{3}k,
\end{align*}
for some positive integer $k$. Then $v>(s+1)^2(t+1)/2$ can be reduced to
\[q_0^3(q_0^3-1)>\frac{(q_0^2-1)^2(q_0^6-q_0^3+1)k^3}{6c}\geq \frac{(q_0^2-1)^2(q_0^6-q_0^3+1)}{18}.\]
This is incorrect for any $q_0$.

Next, we suppose $r\geq 5$. In this case, we have $c=1$ and
\[\frac{t+1}{s+1}=\frac{b}{v}=\frac{3(q_{0}^{2r}-q_{0}^{r}+1)}{q_0^3(q_0^2-1)(q_0^3+1)}=\frac{\frac{q_0^{2r}-q_0^r+1}{q_0^2-q_0+1}}{q_0^3(q_0^2-1)(q_0+1)/3}.\]
Since $\gcd(q_0^{2r}-q_0^r+1, q_{0})=1$, $\gcd(q_0^{2r}-q_0^r+1, q_{0}-1)=1$, $\gcd(q_0^{2r}-q_0^r+1, q_{0}+1)=\gcd(3,q_{0}+1)$ and $q_0^{2r}-q_0^r+1\not\equiv 0\pmod{9}$, then $\gcd\left(\frac{q_0^{2r}-q_0^r+1}{q_0^2-q_0+1},q_0^3(q_0^2-1)(q_0+1)/3\right)=1$. We can assume
\begin{align*}
t+1=\frac{q_0^{2r}-q_0^r+1}{q_0^2-q_0+1}k,\text{~and~} s+1=\frac{q_0^3(q_0^2-1)(q_0+1)}{3}k,
\end{align*}
for some positive integer $k$. Then $v>(s+1)^2(t+1)/2$ is reduced to
\[6q_0^{3r-6}(q_0^{2r}-1)(q_0^r+1)(q_0^2-q_0+1)>(q_0^2-1)^2(q_0+1)^2(q_0^{2r}-q_0^r+1)k^3.\]
Note that $6q_0^{3r-6}(q_0^{2r}-1)(q_0^r+1)^2(q_0^2-q_0+1)<6q_0^{5r-4}(q_0^r+1)^2$, and
$(q_0^2-1)^2(q_0+1)^2(q_0^{3r}+1)k^3>q_0^{3r+6}k^3$. Thus $k^3<6q_0^{2r-10}(q_0^r+1)^2$. This further implies $k<q_0^4$ if $r=5$; $k<q_0^7$ if $r=7$; $k<q_0^{12}$ if $r=11$; and $k<q_0^{15}$ if $r=13$.

Note from \eqref{4.5b} that $q_0^{3r-3}$ divides $b$. Since $b=(t+1)(st+1)$ and $(t+1)_p=(k)_p\leq k$, using the above bounds on $k$, we get $q_0^8\mid st+1$ if $r=5$; $q_0^{11}\mid st+1$ if $r=7$; $q_0^{18}\mid st+1$ if $r=11$; and $q_0^{21}\mid st+1$ if $r=13$.

We first consider the case when $r=13$. Since $k<q_0^{15}$, we assume $k=\sum_{i=0}^{14} a_iq_0^i$, where $a_i\in [0,q_0-1]$. Note that $st+1$ is now expressed in terms of $q_0$ and the $a_i$'s. We can use $q_0^{21}\mid st+1$ repeatedly to deduce some restrictions on the $a_i$'s. At Step $i$, $0\leq i\leq 16$, using $q_0^{i+1}\mid st+1$, we either determine the exact value of $a_i$ or obtain a range for an expression of $a_i$. We verify by computer that for $q_0\leq 236$ there are no integers $s$ and $t$ satisfying both $v=(s+1)(st+1)$ and $b=(t+1)(st+1)$. In the process described above, we suppose $q_0>236$. The results are summarized in Table \ref{tab_cf_r13}. Note that due to technical reasons, for certain $i$'s, we consider $\epsilon(st+1)\pmod{q^i}$ for some integer $\epsilon$ instead of $st+1\pmod{q^i}$.
\begin{table}[h]\aboverulesep=0pt \belowrulesep=0pt
\setlength{\abovecaptionskip}{0cm}
\setlength{\belowcaptionskip}{0cm}
\caption{$T_\a\cong\hat~(q^2-q+1)\!:\!3$, $T_l\cong \PSU(3,q_0). \gcd(3,\frac{q+1}{q_0+1})$, $r=13$, $q_0>236$}
\label{tab_cf_r13}
\centering
\[\scriptsize
\arraycolsep=5pt
\begin{array}{c|ll|ll}
\toprule
\text{Step}  & \text{Modulo condition} & &\text{Claim} &\\
\midrule\midrule
& &\\[-0.3cm]
\text{0.} & st+1\equiv 2- a_0& \pmod{q_0}& a_0=2 &\\
& &\\[-0.3cm]
\text{1.} & st+1\equiv (-2- a_1 ) q_0& \pmod{q_0^2}& a_1=q_0-2& \\
& &\\[-0.3cm]
\text{2.} & 3(st+1)\equiv (3 - 3 a_2) q_0^2 &\pmod{q_0^3}& 3 - 3 a_2=q_0\ld_2 & \ld_2\in[-2,0]\\
& &\\[-0.3cm]
\text{3.} & 3(st+1)\equiv ( \ld_2- 3 a_3  -2) q_0^3 &\pmod{q_0^4}&\ld_2- 3 a_3  -2= q_0\ld_3    & \ld_3\in[-3,-1]\\
& &\\[-0.3cm]
\text{4.} & 3(st+1)\equiv  ( \ld_3- 3 a_4 +2 ) q_0^4& \pmod{q_0^5}& \ld_3- 3 a_4 +2=q_0\ld_4& \ld_4\in[-2,0]\\
& &\\[-0.3cm]
\text{5.} & 3(st+1)\equiv (\ld_4-3a_5) q_0^5 &\pmod{q_0^6}& \ld_4-3a_5=q_0\ld_5    & \ld_5\in[-2,0]\\
& &\\[-0.3cm]
\text{6.}  & 3(st+1)\equiv (\ld_5-3a_6-2) q_0^6 &\pmod{q_0^7}& \ld_5-3a_6-2=q_0\ld_6 & \ld_6\in[-3,-1]\\
& &\\[-0.3cm]
\text{7.}  & 3(st+1)\equiv (\ld_6-3a_7+2) q_0^7 &\pmod{q_0^8}& \ld_6-3a_7+2=q_0\ld_7    & \ld_7\in[-2,0]\\
& &\\[-0.3cm]
\text{8.}  & 27(st+1)\equiv (9\ld_7-27a_8-18) q_0^8 &\pmod{q_0^9}& 9\ld_7-27a_8-18=q_0\ld_8      & \ld_8\in[-27,-1]\\
& &\\[-0.3cm]
\text{9.}  & 27(st+1)\equiv (\ld_8-27a_9+32) q_0^{9} &\pmod{q_0^{10}}& \ld_8-27a_9+32=q_0\ld_9     & \ld_9\in[-26,0]\\
& &\\[-0.3cm]
\text{10.}  & 27(st+1)\equiv ( \ld_9- 27 a_{10} -14 ) q_0^{10} &\pmod{q_0^{11}}& \ld_9- 27 a_{10} -14 =q_0\ld_{10}    & \ld_{10}\in[-27,-1]\\
& &\\[-0.3cm]
\text{11.}  & 27(st+1)\equiv (\ld_{10}- 27 a_{11} -10 ) q_0^{11} &\pmod{q_0^{12}}& \ld_{10}- 27 a_{11} -10 =q_0\ld_{11}  & \ld_{11}\in[-27,-1]\\
& &\\[-0.3cm]
\text{12.}  & 27(st+1)\equiv ( \ld_{11}- 27 a_{12}+6) q_0^{12}& \pmod{q_0^{13}}&  \ld_{11}- 27 a_{12}+6=q_0\ld_{12}  & \ld_{12}\in[-26,0]\\
& &\\[-0.3cm]
\text{13.}  & 27(st+1)\equiv (\ld_{12} - 27 a_{13} + 36) q_0^{13} &\pmod{q_0^{14}}& \ld_{12} - 27 a_{13} + 36=q_0\ld_{13}    & \ld_{13}\in[-26,0]\\
& &\\[-0.3cm]
\text{14.}  & 243(st+1)\equiv (9\ld_{13}- 243 a_{14} ) q_0^{14} &\pmod{q_0^{15}}& 9\ld_{13}- 243 a_{14} =q_0\ld_{14} & \ld_{14}\in[-242,0]\\
& &\\[-0.3cm]
\text{15.}  & 243(st+1)\equiv (\ld_{14}+128) q_0^{15} &\pmod{q_0^{16}}& q_0\mid  \ld_{14}+128 &\\[0.1cm]
\text{16.} & 243(st\!+\!1)\!\equiv\! (\ld_{14} \!+\!128\!  -\!  108 q_0 \! +\!  \ld_{14} q_0)q_0^{15}\!\!\!\!\!\!\!\!\! &\pmod{q_0^{17}}& q_0^2\mid \ld_{14} \!+\!128\!  -\!  108 q_0 \! +\!  \ld_{14} q_0
\\[0.1cm]
\bottomrule
\end{array}
\]
\end{table}

By Step 15, we see $q_0$ divides $\ld_{14}+128$ which lies in $[-114,128]$. Since $q_0>236$, then $\ld_{14}=-128$. So the claim of Step 16 is reduced to $q_0\mid 236$, which contradicts our assumption that $q_0>236$.

For $r\in\{5,7,11\}$, a similar argument will also yield a contradiction. Therefore $T_l$ cannot be isomorphic to $\PSU(3,q_0). \gcd(3,\frac{q+1}{q_0+1})$ and the proof is complete.
\end{proof}

\begin{lemma}\label{d}
The subgroup $T_{\a}$ cannot have type (d), i.e., $T_\a\not\cong\hat~(q+1)^2:\sym_3$.
\end{lemma}
\begin{proof}
Suppose to the contrary that $T_{\a}\cong\hat~(q+1)^2:\sym_3$. By Table \ref{tab_sub}, we have
\[v=q^3(q-1)(q^2-q+1)/6.\]
By the bounds on $[T:T_l]$ in Lemma \ref{lem_boundM1}, we have $q\leq 6730$ if $T_l\cong\hat~3_{+}^{1+2}:\rQ_8.\frac{\gcd(9,q+1)}{3}$; $q\leq 4067$ if $T_l\cong\PSL(2,7)$; and $q\leq 18701$ if $T_l\cong\A_6$. We verify by computer that there are no feasible $(s,t)$ pairs that satisfy $v=(s+1)(st+1)$ and $b=(t+1)(st+1)$ in these cases.

Suppose that $T_l\cong\SO(3,q)$ with $q\geq 7$ odd. By Table \ref{tab_sub} we have $b=q^2(q^3+1)/d$, where $d=\gcd(3,q+1)$. So $\frac{s+1}{t+1}=\frac{v}{b}=\frac{q(q-1)d}{6(q+1)}$.  Let $\theta=\gcd(q(q-1)d/6,q+1)=\gcd((q-1)/2,2)$. Then there exists a positive integer $k$ such that
\begin{equation}\label{dest}
s+1=\frac{q(q-1)d}{6\theta}k,\text{ and }t+1=\frac{q+1}{\theta}k.
\end{equation}
By Lemma \ref{newcond}\ref{cond1} we have
$
q^3(q-1)(q^2-q+1)/6=v>(s+1)^2(t+1)/2\geq \frac{q^2(q-1)^2(q+1)k^3}{72\theta^3}$, which implies $12q(q^2-q+1)\theta^3\geq (q^2-1)k^3$. If $k\geq q$, then $96(q^2-q+1)>(q^2-1)q^2$ which holds only if $q=9$. But there are no integers $s$ and $t$ satisfying $v=(s+1)(st+1)$ and $b=(t+1)(st+1)$ for $q=9$. Thus $k<q$. Recall that $b=q^2(q^3+1)/d$. This implies that $q^2\mid b$. Since $b=(t+1)(st+1)$ and $(t+1)_p=(k)_p\leq k<q$, then $q\mid st+1$. By \eqref{dest}, this is reduced to $q\mid 6\theta(k-2\theta)$. Note that $q$ is odd and $\theta=1$ or $2$. Then $q\mid 3(k-2\theta)$. Since $k<q$, we have $3(k-2\theta)\in\{0,q,2q\}$ and so $k\in\{2\theta,2\theta+\frac{q}{3},2\theta+\frac{2q}{3}\}$. We check by computer that there are no prime power $q$ satisfying $b=(t+1)(st+1)$ no matter which expression $k$ takes. Therefore $T_l$ cannot be isomorphic to $\SO(3,q)$. 

Suppose that $T_l\cong\PSU(3,q_0). \gcd(3,\frac{q+1}{q_0+1})$, where $q=q_0^r$, $r$ odd and prime. By Table \ref{tab_sub}, we have $b=\frac{q^3(q^2-1)(q^3+1)}{c\cdot q_0^3(q_0^2-1)(q_0^3+1)}$. By Lemma \ref{newcond}\ref{cond2}, we have
\begin{align*}
v>\frac{1}{2}\left(\frac{b}{v}\right)^4=\frac{1}{2}\left(\frac{6}{c}\cdot\frac{(q_0^r+1)^2}{q_0^3(q_0^2-1)(q_0^3+1)}\right)^4
> \frac{1}{2}\left(\frac{2(q_0^r+1)^2}{q_0^8}\right)^4>  8q_0^{8r-32}.
\end{align*}
Note that $v=q_0^{3r}(q_0^{r}-1)(q_0^{2r}-q_0^r+1)/6<q_0^{6r}/6$. Thus $q_0^{32-2r}>48$ which holds only if $r\in\{3,5,7,11,13\}$.

We first consider the case when $r=3$. We have $\frac{t+1}{s+1}=\frac{b}{v}=\frac{q_0^2-q_0+1}{cq_0^3(q_0-1)/6}$. Since $\gcd(q_0^2-q_0+1,cq_0^3(q_0-1)/6)=1$, there exists an integer $k$ such that $t+1=(q_0^2-q_0+1)k$ and $s+1=q_0^3(q_0-1)k\cdot \frac{c}{6}$. If $k\geq q_0^3$, then we have
\begin{align*}
q_0^{18}/6>&q_0^9(q_0^3-1)(q_0^6-q_0^3+1)/6=v>(s+1)^2(t+1)/2\\
=&\left(q_0^3(q_0-1)\cdot \frac{c}{6}\right)^2(q_0^2-q_0+1)k^3/2\geq q_0^{15}(q_0-1)^2(q_0^2-q_0+1)/72,
\end{align*}
which implies $12q_0^3>(q_0-1)^2(q_0^2-q_0+1)$, which holds only if $q_0\leq 13$. However, there are no feasible $(s,t)$ pairs that satisfy $v=(s+1)(st+1)$ and $b=(t+1)(st+1)$ for $r=3$ and $q_0\leq 13$. Thus $k<q_0^3$, so we can assume $k=a_0+a_1q_0+a_2q_0^2$ with $0\leq a_0,a_1,a_2<q$. Recall that $b=\frac{q^3(q^2-1)(q^3+1)}{c\cdot q_0^3(q_0^2-1)(q_0^3+1)}$. Then $q_0^{6}\mid b$. Since $b=(t+1)(st+1)$ and $(t+1)_p=(k)_p\leq k<q_0^3$, we deduce that $q_0^{3}\mid st+1$. Computing $st+1$ modulo $q_0^i$, $1\leq i\leq 3$, we obtain $a_0=a_1=2$, and $a_2\in\{0,2q_0,3q_0\}$. Now, both $s$ and $t$ are expressions of $q_0$. We check by computer that $v=(s+1)(st+1)$ has no prime power solution for $q_0$.

We then consider the case when $5\leq r\leq13$. In this case, $c$ is always 1. We have  $\frac{t+1}{s+1}=\frac{b}{v}=\frac{6(q_{0}^{r}+1)^{2}}{q_0^3(q_{0}^{2}-1)(q_{0}^{3}+1)}$.  Since $\gcd(q_{0}^{r}+1, q_{0}^{2}-1)=q_0+1$ and $\gcd(q_{0}^{r}+1, q_{0}^{3}+1)=q_0+1$, there exists an integer $k$ such that
\[t+1=\left(\frac{q_0^r+1}{q_0+1}\right)^2k, \text{ and }s+1=\frac{q_0^3(q_{0}^{2}-1)(q_{0}^{3}+1)k}{(q_0+1)^2}=q_0^3(q_0-1)(q_0^2-q_0+1)k/6.\]
By Lemma \ref{newcond}\ref{cond1}, we have
\begin{align*}
v&>(s+1)^2(t+1)/2=\frac{1}{72}(q_0^3(q_0-1)(q_0^2-q_0+1))^2\left(\frac{q_0^r+1}{q_0+1}\right)
^2k^3  \\
&>\frac{1}{72} \left(q_0^3\cdot \frac{q_0}{2}\cdot \frac{q_0^2}{2}\right)^2  \cdot \frac{q_0^{2r}}{4q_0^2} k^3=\frac{q_0^{2r+10}k^3}{72\cdot 64}.
\end{align*}
Note that $v=q_0^{3r}(q_0^r-1)(q_0^{2r}-q_0^r+1)/6<q_0^{6r}/6$. It follows that $k^3<768q_0^{4r-10}$.

\begin{table}[h]\aboverulesep=0pt \belowrulesep=0pt
\setlength{\abovecaptionskip}{0cm}
\setlength{\belowcaptionskip}{0cm}
\caption{$T_\a\cong\hat~(q+1)^2:\sym_3$, $T_l\cong \PSU(3,q_0). \gcd(3,\frac{q+1}{q_0+1})$, $r=13$, $q_0>29988$}
\label{tab_df_r13}
\centering
\[\scriptsize
\arraycolsep=5pt
\hspace{-.4cm}
\begin{array}{c|ll|ll}
\toprule
&&&&\\[-0.3cm]
\text{Step} & \text{Modulo condition} & & \text{Claim} &\\
\midrule\midrule
& &\\[-0.3cm]
\text{0.} & st+1\equiv 2- a_0 &\hspace{-0.4cm}\pmod{q_0}& a_0=2 &\\
& &\\[-0.3cm]
\text{1.} & st+1\equiv (4- a_1 ) q_0& \hspace{-0.4cm}\pmod{q_0^2}& a_1=4& \\
& &\\[-0.3cm]
\text{2.} & 3(st+1)\equiv (6 - 3 a_2) q_0^2& \hspace{-0.4cm}\pmod{q_0^3}& 6-3a_2=q_0\ld_2 & \ld_2\in[\text{-}2,0]\\
& &\\[-0.3cm]
\text{3.} & 3(st+1)\equiv ( \ld_2- 3 a_3-1 ) q_0^3 &\hspace{-0.4cm}\pmod{q_0^4}& \ld_2 - 3 a_3 -1 =q_0\ld_3       & \ld_3\in[\text{-}3,\text{-}1]\\
& &\\[-0.3cm]
\text{4.} & 3(st+1)\equiv  (\ld_3 - 3 a_4-2 ) q_0^4 &\hspace{-0.4cm}\pmod{q_0^5}& \ld_3 - 3 a_4-2=q_0\ld_4      & \ld_4\in[\text{-}3,\text{-}1]\\
& &\\[-0.3cm]
\text{5.} & 6(st+1)\equiv (2 \ld_4-6 a_5) q_0^5 &\hspace{-0.4cm}\pmod{q_0^6}&2 \ld_4-6 a_5= q_0\ld_5          & \ld_5\in[\text{-}6,\text{-}1]\\
& &\\[-0.3cm]
\text{6.}  & 6(st+1)\equiv (\ld_5 - 6 a_6 +3 ) q_0^6 &\hspace{-0.4cm}\pmod{q_0^7}& \ld_5 - 6 a_6 +3=q_0\ld_6& \ld_6\in[\text{-}5,0]\\
& &\\[-0.3cm]
\text{7.}  & 6(st+1)\equiv (\ld_6-6 a_7) q_0^7 &\hspace{-0.4cm}\pmod{q_0^8}& \ld_6-6 a_7=q_0\ld_7 & \ld_7\in[\text{-}5,0]\\
& &\\[-0.3cm]
\text{8.}  & 108(st+1)\equiv (18\ld_7- 108 a_8 -18 ) q_0^8 &\hspace{-0.4cm}\pmod{q_0^9}& 18\ld_7- 108 a_8 -18=q_0\ld_8&  \ld_8\in[\text{-}108,\text{-}1]\\
& &\\[-0.3cm]
\text{9.}  & 108(st+1)\equiv ( \ld_8 - 108 a_9 +25) q_0^{9} &\hspace{-0.4cm}\pmod{q_0^{10}}& \ld_8 - 108 a_9 +25=q_0\ld_9   & \ld_9\in[\text{-}107,0]\\
& &\\[-0.3cm]
\text{10.}  & 108(st+1)\equiv ( \ld_9- 108 a_{10} +68 ) q_0^{10} &\hspace{-0.4cm}\pmod{q_0^{11}}& \ld_9- 108 a_{10} +68 =q_0\ld_{10} & \ld_{10}\in[\text{-}107,0]\\
& &\\[-0.3cm]
\text{11.}  & 216(st+1)\equiv ( 2\ld_{10}- 216 a_{11} +44 ) q_0^{11} &\hspace{-0.4cm}\pmod{q_0^{12}}& 2\ld_{10}- 216 a_{11} +44=q_0\ld_{11} & \ld_{11}\in[\text{-}215,0]\\
& &\\[-0.3cm]
\text{12.}  & 216(st+1)\equiv (\ld_{11}- 216 a_{12}-135 ) q_0^{12} &\hspace{-0.4cm}\pmod{q_0^{13}}& \ld_{11}- 216 a_{12}-135=q_0\ld_{12} & \ld_{12}\in[\text{-}216,\text{-}1]\\
& &\\[-0.3cm]
\text{13.}  & 216(st+1)\equiv ( \ld_{12}- 216 a_{13} -966 ) q_0^{13} &\hspace{-0.4cm}\pmod{q_0^{14}}& \ld_{12}- 216 a_{13} -966=q_0\ld_{13} & \ld_{13}\in[\text{-}216,\text{-}1]\\
& &\\[-0.3cm]
\text{14.}  & 3888(st\!+\!1)\!\equiv \!(18\ld_{13} \!-  \!3888 a_{14}  \!- \!29970) q_0^{14}&\hspace{-0.4cm}\pmod{q_0^{15}}&  \!18\ld_{13} \!- \!3888 a_{14} \! - \!29970\!=\!q_0\ld_{14}& \ld_{14}\in[\text{-}3888,\text{-}1]\\
& &\\[-0.3cm]
\text{15.}  & 3888(st+1)\equiv (\ld_{14}-12977 ) q_0^{15} &\hspace{-0.4cm}\pmod{q_0^{16}}& q_0\mid \ld_{14}-12977 &\\[0.1cm]
\bottomrule
\end{array}
\]
\end{table}
Suppose $r=13$. We check by computer that there are no integers $s$ and $t$ satisfying both $v=(s+1)(st+1)$ and $b=(t+1)(st+1)$ if $q_0\leq 29988$. For the rest of the proof, we assume $q_0>29988$. By $k^3<768q_0^{4r-10}$, we have $k<q_0^{15}$. Recall that $b=\frac{q^3(q^2-1)(q^3+1)}{c\cdot q_0^3(q_0^2-1)(q_0^3+1)}=\frac{q_0^{36}(q_0^{26}-1)(q_0^{39}+1)}{(q_0^2-1)(q_0^3+1)}$. Then $q_0^{36}\mid b$. Since $b=(t+1)(st+1)$ and $(t+1)_p=k_p\leq k<q_0^{15}$, we have $q_0^{21}\mid st+1$. We repeat the process as in Table \ref{tab_cf_r13}. The results are summarized in Table \ref{tab_df_r13}. By the last step, we have $q_0$ divides $\ld_{14}-12977$ which lies in $[-16865, -12978]$. Thus $q_0\leq 16865$, but this is contrary to our assumption that $q_0>29988$.

For $r\in\{5,7,11\}$, a similar argument will also yield a contradiction. Therefore $T_l$ cannot be isomorphic to $\PSU(3,q_0). \gcd(3,\frac{q+1}{q_0+1})$ and the proof is complete.
\end{proof}

\begin{lemma}\label{e}
The subgroup $T_{\a}$ cannot have type (e), i.e., $T_\a\not\cong\SO(3,q)$ with $q\geq7$ odd.
\end{lemma}
\begin{proof}
Suppose to the contrary that $T_{\a}\cong\SO(3,q)$ with $q\geq7$ odd. By Table \ref{tab_sub}, we have
\begin{equation}\label{ev}
v=q^2(q^3+1)/d,
\end{equation}
where $d=\gcd(3,q+1)$. By the bounds on $[T:T_l]$ in Lemma \ref{lem_boundM1}, we have $q\leq17$ if $T_l\cong\hat~3_{+}^{1+2}:\rQ_8.\frac{\gcd(9,q+1)}{3}$; $q\leq 19$ if $T_l\cong\PSL(2,7)$; and $q=11$ if $T_l\cong\A_6$. We verify by computer that there are no feasible $(s,t)$ pairs that satisfy $v=(s+1)(st+1)$ and $b=(t+1)(st+1)$ in these cases.

Suppose that $T_l\cong\PSU(3,q_0). \gcd(3,\frac{q+1}{q_0+1})$, where $q=q_0^r$, $r$ odd and prime. We see from Table \ref{tab_sub} that 
\begin{equation}\label{fb}
b=\frac{q^3(q^2-1)(q^3+1)}{c\cdot q_0^3(q_0^2-1)(q_0^3+1)},
\end{equation}
where $c=\gcd(3,\frac{q+1}{q_0+1})$. Note that $c\leq\gcd(q+1,3)=d$. By Lemma \ref{newcond}\ref{cond2}, we have
\begin{align*}
v&>\frac{1}{2}\left(\frac{b}{v}\right)^4=\frac{1}{2}\left(\frac{d}{c}\cdot \frac{q(q^2-1)}{q_0^3(q_0^2-1)(q_0^3+1) }\right)^4>\frac{q^4(q^2-1)^4}{2q_0^{32} }>\frac{q^{12}}{4q_0^{32}}.
\end{align*}
Since $v=q^2(q^3+1)/d<2q^5$, we deduce that $8>q_0^{7r-32}$. This implies that $r=3$. In this case, we have $d=c$, and $\frac{t+1}{s+1}=\frac{b}{v}=\frac{q_0^2+q_0+1}{q_0+1}$. Since $\gcd(q_0^2+q_0+1,q_0+1)=1$, there exists an integer $k$ such that
\[t+1=(q_0^2+q_0+1)k, \text{ and }s+1=(q_0+1)k.\]
If $k\geq q_0^4$, then
$
(s+1)^2(t+1)/2\geq(q_0+1)^2(q_0^2+q_0+1)q_{0}^{12}/2
> q_0^6(q_0^9+1)\geq v$,
contradicting Lemma \ref{newcond}\ref{cond1}. If $k<q_0^2$, then
$
q_0^6(q_0^9+1)/3\leq v<(s+1)^2(t+1)<(q_0+1)^2(q_0^2+q_0+1)q_{0}^6$,
a contradiction. Thus \[q_0^2\leq k< q_0^4.\] Since $q_0^6\mid v=(s+1)(st+1)$ and $(s+1)_p=(k)_p\leq k<q_0^4$, it follows that
$q_{0}^{2}\mid st+1=(q_0k+k-1)(q_0^2k+q_0k+k-1)+1$. Thus $q_0\mid(k-1)^2+1$. We claim that $\gcd(s,t)=1$. Indeed, noting that $\gcd(s,t)=\gcd((q_0+1)k-1,k+q_0-1)=\gcd((q_0+1)k-1,q_0^2)$, if $\gcd(s,t)\neq1$, then $p$ divides $(q_0+1)k-1$, and so $p$ divides $k-1$, a contradiction. Thus $\gcd(s,t)=1$. Then we have
\begin{equation}\label{efgcd1}
\gcd(s,s+t)=\gcd(t,s+t)=1.
\end{equation} By \eqref{ev} and \eqref{fb}, we see that $v$ is even and $b$ is odd. Thus $s,k$ are odd and $t$ is even. It follows that
\begin{equation}\label{efgcd2}
\gcd(s+t,k)=\gcd((q_0^2+2q_0+2)k-2,k)=\gcd(2,k)=1.
\end{equation}
Note that $st(s+1)(t+1)=st(q_0^2+q_0+1)(q_0+1)k^2$. Then by \eqref{efgcd1} and \eqref{efgcd2}, the divisibility condition in Lemma \ref{para}\ref{divcond} is reduced to $s+t\mid (q_0^2+q_0+1)(q_0+1)$. In particular \[s+t=(q_0^2+2q_0+2)k-2<(q_0^2+q_0+1)(q_0+1).\]
On the other hand, since $k\geq q_0^2$, we have $s+t\geq(q_0^2+2q_0+2)q_0^2-2> (q_0^2+q_0+1)(q_0+1)$, a contradiction. This completes the proof.
\end{proof}

\begin{lemma}\label{f}
The subgroup $T_{\a}$ cannot have (f), i.e., $T_\alpha\ncong\PSU(3,q_0). \gcd(3,\frac{q+1}{q_0+1})$,  $q=q_0^r$, $r$ odd and prime.
\end{lemma}
\begin{proof}
This follows directly from the restrictions on $q$. Since $T_\alpha$ has type (f), then $q$ is not a prime. However, since $T_l$ has type (g), (h) or (i), then $q$ is a prime, a contradiction.
\end{proof}

\begin{lemma}\label{gh}
The subgroup $T_{\a}$ cannot have type (g) or (h), i.e, $T_\a\not\cong\hat~3_{+}^{1+2}:\rQ_8.\frac{\gcd(9,q+1)}{3}$ with $p=q\equiv 2\pmod 3$, $q\geq 11$, and $T_\a\not\cong\PSL(2,7)$ with $q\neq 5$ and $p=q\equiv 3,5,6\pmod{7}$.
\end{lemma}
\begin{proof}
The proofs with respect to the two types are similar. We only prove that $T_\a$ cannot have type (g) here.

Suppose to the contrary that $T_{\a}=\hat~3_{+}^{1+2}:\rQ_8.\frac{\gcd(9,q+1)}{3}$ with $p=q\equiv 2\pmod 3$, $q\geq 11$. Then $v=\frac{q^3(q^2-1)(q^3+1)}{72\cdot \gcd(9,q+1)}$ by Table \ref{tab_sub}. Note that $v$ is even, so $s$ is odd.

We first assume $T_l\cong\PSL(2,7)$ with $q\neq 5$, $p=q\equiv3,5,6\pmod 7$. It follows form Table \ref{tab_sub} that $b=\frac{q^3(q^2-1)(q^3+1)}{168\cdot\gcd(3,q+1)}$. Then $\frac{s+1}{t+1}=\frac{7}{\gcd(9,q+1)}$. So there exists a positive integer $k$ such that
\[s+1=7k, \text{ and }t+1=\begin{cases}9k,& \text{if }9\mid q+1,\\ 3k, &\text{if }9\nmid q+1.\end{cases}\]
Since $s$ is odd, then $k$ must be even. If $9\mid q+1$, then the divisibility condition in Lemma \ref{para}\ref{divcond} is reduced to $16k-2\mid63k^2(7k-1)(9k-1)$, which simplifies to $8k-1\mid 63$. Thus $k=8$, $s=55$ and $t=71$. But there is no prime $q=p$ that satisfies $v=(s+1)(st+1)$ and $b=(t+1)(st+1)$. If $9\nmid q+1$, then the divisibility condition in Lemma \ref{para}\ref{divcond} is reduced to $10k-2\mid 21k^2(7k-1)(3k-1)$. This simplifies to $5k-1\mid 42$,
which holds for no positive even integers. Thus $T_l\not\cong\PSL(2,7)$.

Next, we suppose $T_l\cong\A_6$ with $p=q\equiv 11,14\pmod {15}$. It follows from Table \ref{tab_sub} that  $\frac{s+1}{t+1}=\frac{v}{b}=\frac{15}{\gcd(9,q+1)}$. Note that $\gcd(9,q+1)=3$ or $9$ accoring as $9$ divides $q+1$ or not. So there exists an integer $k$ such that \[s+1=5k, \text{ and }t+1=\begin{cases}3k,&\text{if }9\mid q+1, \\k, &\text{if }9\nmid q+1.\end{cases}\]
 Since $s$ is odd, $k$ is even. If $9\mid q+1$, then the divisibility condition in Lemma \ref{para}\ref{divcond} is reduced to $8k-2\mid 15k^2(5k-1)(3k-1)$, which simplifies to $4k-1\mid 15$. Thus $k=4$, $s=19$ and $t=11$. But there is no prime $q=p$ that satisfies $v=(s+1)(st+1)$ and $b=(t+1)(st+1)$. If $9\nmid q+1$, then the divisibility condition in Lemma \ref{para}\ref{divcond} is reduced to $6k-2\mid 5k^2(5k-1)(k-1)$. It simplifies to $3k-1\mid 10$. Thus $k=2$, $s=9$ and $t=1$. This contradicts the assumption that $t>1$. Thus $T_l\not\cong\A_6$. The proof is complete.
\end{proof}

\subsection{The case when $T_\a$ and $T_l$ have the same type}

\begin{lemma}\label{bb}
The subgroup $T_{\a}$ cannot have type (b), i.e., $T_\a\not\cong\hat~\GU(2,q)$.
\end{lemma}
\begin{proof}
Suppose to the contrary that $T_\alpha\cong\hat~\GU(2,q)$. By \eqref{point} and Table \ref{tab_sub}, we have
\begin{equation}\label{bb1}
v=(s+1)(s^2+1)=q^2(q^2-q+1).
\end{equation}

We first consider the case when $q$ is odd. The above equation implies that $s$ is even. So $\gcd(s+1,s^2+1)=\gcd(s+1,2)=1$. It follows that $(v)_p=(s+1)_p=q^2$ or $(v)_p=(s^2+1)_p=q^2$. If $q^2\mid s+1$, then $(s+1)(s^2+1)\geq q^2(q^4-2q^2+2)>q^2(q^2-q+1)=v$, a contradiction. Hence $q^2\mid s^2+1$. We assume that there exists some integer $x$ such that
\begin{align}
\label{bb2}s^2+1=q^2x,\\
\label{bb3}q^2-q+1=(s+1)x.
\end{align}
By \eqref{bb1} and \eqref{bb2} we have
\begin{equation}\label{bb4}
x^3q^6=(s^2+1)^3<(s+1)^2(s^2+1)^2=q^4(q^2-q+1)^{2}<q^8.
\end{equation}
Thus $x^{3}<q^{2}$. Solving $x$ in \eqref{bb3} and substituting it in \eqref{bb2}, we get $(q^2-q+1-x)^2+x^2=q^2x^3$ which implies $q^2\mid 2q(x-1)+(x-1)^2+x^2$. In particular, we have $q^2<2q(x-1)+(x-1)^2+x^2$. This further implies \[q^6\leq (2q(x-1)+(x-1)^2+x^2)^3<(2qx+2x^2)^3=8x^3(q+x)^3<8q^2(q+q)^3=64q^5.\]
Thus $q<64$. It is verified by computer that \eqref{bb1} has no integer solutions for $s$ when $q<64$.

We then consider the case when $q$ is even. By \eqref{bb1}, we see that $s$ is odd. Thus $(s^2+1)_2=2$ and $(s+1)_2=q^2/2$. This implies $(s+1)(s^2+1)\geq q^2(q^4-4q^2+8)/8>v$, a contradiciton. This completes the proof.
\end{proof}

\begin{lemma}\label{cc}
The subgroup $T_\alpha$ cannot have type (c), (d) or (e), i.e., $T_\a\not\cong\hat~(q^2-q+1):3$ with $q\neq 3,5$, $T_\a\not\cong\hat~(q+1)^2:\sym_3$ with $q\neq 5$, and $T_\a\not\cong\SO(3,q)$ with $q\geq 7$, $q$ odd.
\end{lemma}
\begin{proof}
The proofs of these cases are similar, so we only prove that $T_\a$ cannot have type (c) here. Suppose to the contrary that $T_{\a}\cong\hat~(q^2-q+1):3$ with $q\neq 3,5$. By Table \ref{tab_sub}, we have
\begin{equation}\label{cc1}
3v=q^3(q^2-1)(q+1)=3(s+1)(s^2+1).
\end{equation}
It is verified by computer that \eqref{cc1} has no integer solutions for $s$ when $q<104$. Next, we assume $q\geq 104$ which consequently implies that $s\geq 7523$. 

We first consider the case when $q$ is odd and $q\not\equiv 0\pmod{3}$. Then  $s$ is odd and $(s+1,s^2+1)=2$. This means $q^{3}\mid s+1$ or $q^{3}\mid s^{2}+1$. If $q^{3}\mid s+1$, then $3(s+1)(s^2+1)>3q^3(q^6-2q^3+2)>q^3(q^2-1)(q+1)$. Thus $q^{3}\mid s^{2}+1$ and we assume that there exist some integers $x$ such that
\begin{align}
\label{cc2}s^2+1=q^3x,\\
\label{cc3}(q^2-1)(q+1)=3(s+1)x.
\end{align}
Note that \eqref{cc1} implies
 \[q^6<q^3(q^2-1)(q+1)=3(s+1)(s^2+1)<4s^3,\]
\[2q^6>q^3(q^2-1)(q+1)=3(s+1)(s^2+1)>3s^3.\]
Thus we have $q^6/4<s^3<2q^6/3$. Combining \eqref{cc2} and this bound,  we get
\[(q^3x)^3>s^6>q^{12}/16,\]
\[(q^3x)^3<2s^6<8q^{12}/9.\]
It follows that $q/3<10q/27<q/\sqrt[3]{16}<x<q$. Solving for $x$ in \eqref{cc3} and substituting it in \eqref{cc2}, we get
\[((q^2-1)(q+1)-3x)^2+9x^2=9x^3q^3,\]
which implies
\[9x^2+(q+1+3x)^2-2q^2(3x+1)\equiv 0\pmod{q^3}.\]
By the fact that $q/3<x<q$, we have
\[-6q^3<9x^2+(q+1+3x)^2-2q^2(3x+1)<-q^3.\]
Moreover, since both $q$ and $s$ are odd, then \eqref{cc2} implies that $x$ is even and so $9x^2+(q+1+3x)^2-2q^2(3x+1)\equiv2\pmod{4}$. Thus $9x^2+(q+1+3x)^2-2q^2(3x+1)=-2q^3$.  We rewrite this as $2q^2(3x+1)-2q^3=9x^2+(q+1+3x)^2$. Note that
\begin{align*}
2q^2(3x+1)-2q^3>2q^2(3\cdot10q/27+1)-2q^3=&2q^3/9+2q^2\\
>&9q^2+(4q+1)^2>9x^2+(q+1+3x)^2
\end{align*}
for $q\geq 104$. We conclude that $9x^2+(q+1+3x)^2-2q^2(3x+1)=-2q^3$ has no integer solutions. Therefore \eqref{cc1} has no integer solutions.

We then consider the case $q\equiv 0\pmod{3}$. By a similar argument, we get $\frac{q^3}{3}\mid s^2+1$. Thus $3$ divides $s^2+1$ which is impossible. Finally, we assume that $q$ is even. In this case, we see that $s$ is odd, so $q^3/2$ strictly divides $s+1$. As a result, we have $3(s+1)(s^2+1)\geq 3q^3(q^6-4q^3+16)/8>q^3(q^2-1)(q+1)$, a contradiction. This completes the proof.
\end{proof}

\begin{lemma}\label{-1square}
Suppose $p$ is odd and $T_\alpha$ has type (f), (g), (h) or (i). Then $p\equiv 1\pmod{4}$.
\end{lemma}
\begin{proof}
Noting that $\gcd(s+1,s^2+1)=\gcd(s+1,2)\leq 2$, we have $(v)_p=(s+1)_p$ or $(v)_p=(s^2+1)_p$.

Suppose that $T_\alpha$ has type (f) with $p$ odd, then by Table \ref{tab_sub}, we have $v=\frac{q^3(q^2-1)(q^3+1)}{c\cdot q_0^3(q_0^2-1)(q_0^3+1)}$, and so $(v)_p=q_0^{3r-3}$. If $(s+1)_p=q_0^{3r-3}$, then $(s+1)(s^2+1)\geq q_0^{3r-3}(q_0^{6r-6}-2q_0^{3r-3}+2)>v$, a contradiction. Thus $(s^2+1)_p=q_0^{3r-3}$, and so $-1\equiv s^2\pmod{p}$. Therefore $p\equiv 1\pmod{4}$.

If $T_\alpha$ has type (g), (h) or (i), then we have $p\neq 2,5,7$. We claim that $p\neq 3$. In fact, $p=3$ can appear only if $T_\a$ has type (h) in which case $v=(s+1)(s^2+1)$ has no integer solutions for $s$. So $p\neq 2,3,5,7$. This implies $(v)_p=q^3$. By a similar argument as in the first part, we obtain $(s^2+1)_p=q^3$ and so $p\equiv 1\pmod{4}$.
\end{proof}

\begin{lemma}\label{ffqeven}
Suppose $T_\a\cong\PSU(3,q_0).\gcd(3,\frac{q+1}{q_0+1})$, where $q=q_0^r$ even and $r$ odd prime. Let $g\in T_\a$ be of order 4. Then $g^T\cap T_\a=g^{T_\a}$. 
\end{lemma}
\begin{proof}
The statement is proved in the following three cases.

(1) Assume $q_0\equiv 1\pmod{3}$. 

Then $\gcd(3,\frac{q+1}{q_0+1})=1$. By Table \ref{tab_ord4} we see that $T_\a\cong\PSU(3,q_0)$ has a unique conjugacy class of elements of order 4. It is clear that $g^T\cap T_\a=g^{T_\a}$. 

(2) Assume $q_0\equiv 2 \pmod{3}$ and $r=3$. 

In this case, we have $T_\a=K.3$ where $K\cong\PSU(3,q_0)$. By Table \ref{tab_ord4}, we see that $K$ has 3 conjugacy classes of elements of order 4. Let $\Omega$ be the set of elements of $T_\a$ of order 4. Then $\Omega$ is contained in $K$. In other words, the elements of $T_\a$ of order 4 are exactly the elements of $K$ of order 4. Consider the actions of $K$ and $T_\a$ on $\Omega$ by conjugation, respectively. If $T_\a$ has 3 orbits, then the $T_\a$-orbits are precisely the $K$-orbits. Let $h\in\Omega$. We have $\frac{|K|}{|C_K(h)|}=\frac{|T_\a|}{|C_{T_\a}(h)|}$. Since $T_\a=K.3$, then $|C_{T_\a}(h)|=3\cdot|C_K(h)|$. Note that $C_{T_\a}(h)\leq C_T(h)$, and $|C_T(h)|=4^f$ (see Table \ref{tab_ord4}). Thus $|C_{T_\a}(h)|$ cannot have 3 as a divisor, so $T_\a$ cannot have 3 orbits on $\Omega$. Similarly, we can prove that $T_\a$ cannot have 2 orbits on $\Omega$. Therefore $T_\a$ is transitive on $\Omega$, i.e., all the elements of $T_\a$ of order 4 are conjugate. The statement follows.

(3) Assume $q_0\equiv 2 \pmod{3}$ and $r>3$. 

Since $r>3$, then $\gcd(3,\frac{q+1}{q_0+1})=1$. Thus $T_\a\cong\PSU(3,q_0)$. Following Table \ref{tab_ord4}, we define 
\[X=\{A\in\SL(3,q^2)\mid AW\overline{A}^\top=W\}\text{ and }X_0=\{A\in\SL(3,q_0^2)\mid AW\widetilde{A}^\top=W\},\]
where $W={\footnotesize\begin{pmatrix}0&0 &1\\ 0&1&1\\ 1&1 &0\end{pmatrix}}$ and $\widetilde{a}=a^{q_0}$. Then $\hat~X\cong T$ and $\hat~X_0\cong T_\a$. For $\ell\in\{1,2,3\}$, define $g_\ell$ to be the image of ${\footnotesize \begin{pmatrix}1&0&0\\ z_\ell &1 &0 \\ 0 & \widetilde{z_\ell} &1\end{pmatrix}}$ in $\hat~X_0$ where $z_\ell=\w_0^{(q_0-1)\ell}+1$ and $\w_0=\w^{\frac{q^2-1}{q_0^2-1}}$. Then $g_1,g_2$ and $g_3$ are pairwise non-conjugate elements of $\hat~X_0$. We claim that they are also pairwise non-conjugate elements of $\hat~X$. First of all, note that $z_\ell=\w_0^{(q_0-1)\ell}+1=(\w^{q-1})^{\frac{q+1}{q_0+1}\cdot\ell}+1=\rho^{\frac{q+1}{q_0+1}\cdot\ell}+1$. Furthermore, since $\w_0^q=\w_0^{q_0^r}=\w_0^{q_0}$, then $\widetilde{z_\ell}=\overline{z_\ell}$. Thus $g_1,g_2,g_3$ are indeed contained in $\hat~X$. Since $\gcd(3,\frac{q+1}{q_0+1})=1$, then $\frac{q+1}{q_0+1}i\not\equiv\frac{q+1}{q_0+1}j\pmod{3}$ whenever $i\not\equiv j\pmod{3}$. According to Table \ref{tab_ord4}, this shows that $g_1,g_2,g_3$ are pairwise non-conjugate in $\hat~X$. Moreover, two non-conjugate elements of $\hat~X_0$ are not conjugate in $\hat~X$. By \cite[Table 8.5]{lowdim}, we see that $T$ has a unique conjugacy class of $\PSU(3,q_0)$. Since $T\cong \hat~X$, we conclude that two non-conjugate elements of $T_\a$ of order 4 are not conjugate in $T$. This implies $g^T\cap T_\a=g^{T_\a}$. 
\end{proof}

\begin{lemma}\label{ff}
The subgroup $T_{\a}$ cannot have type (f), i.e., $T_\a\not\cong\PSU(3,q_0).\gcd(3,\frac{q+1}{q_0+1})$, where $q=q_0^r$ and $r$ odd prime.
\end{lemma}
\begin{proof}
Suppose to the contrary that $T_{\a}=K.c$, where $K \cong \PSU(3,q_0)$ and $c=\gcd(3,\frac{q+1}{q_0+1})$.

We first consider the case when $q$ is odd.  

By Lemma \ref{-1square} we have $q\equiv 1 \pmod 4$. Note that $c=1$ or $3$. So all the elements of $K.c$ of order $4$ lie in $K$. By Table \ref{tab_ord4}, we see that there is a unique conjugacy class of elements of order $4$ in $K$, $T_\alpha$, $T_l$ and $T$. Let $g\in T_\a$ and $h\in T_l$ be of order 4. Then there exists $x\in T$ such that $g^x=h$. It follows that $h\in T_{\a}^x\cap T_l=T_{\b}\cap T_l$, where $\b=\a^x$. Furthermore, we have $h^T\cap T_\beta=h^{T_\beta}$ and $h^T\cap T_l=h^{T_l}$. 
 
Note that $h^{T_\beta}=h^{K^x}$, so we have $|C_{T_\beta}(h)|/|C_{K^x}(h)|=|T_\beta|/ |K^x|=c$. Then by Lemma \ref{order4}\ref{ord4.1}, we have
\[\frac{|C_{T}(h)|}{|C_{T_\beta}(h)|}=\frac{|C_T(h)|}{c\cdot |C_{K^x}(h)|}=\frac{q^2-1}{c(q_0^2-1)}>4.\] Similarly, we have $|C_T(h)|/|C_{T_l}(h)|>4$. By Corollary \ref{centralizer}, we have $C_T(h)$ is nonabelian which is a contradiction since $C_T(h)$ is cyclic as shown in Lemma \ref{order4}\ref{ord4.1}.

Next, we consider the case when $q$ is even and $s=t$. 

By Table \ref{tab_sub}, we have $v=[T:T_{\a}]=\frac{q^3(q^2-1)(q^3+1)}{c\cdot q_0^3(q_0^2-1)(q_0^3+1)}$ which is even. On the other hand, since $v=(s+1)(s^2+1)$, then $s$ is odd and $(s^2+1)_2=2$. Thus $s+1\geq q_0^{3r-3}/2$. If $q_0\geq 4$, we have
\begin{align*}
(s\!+\!1)(s^2\!+\!1)\geq q_0^{3r-3}\Big(\frac{q_0^{6r-6}}{8}\!-\!\frac{q_0^{3r-3}}{2}\!+\!1\!\Big)\!>q_0^{3r-3}\Big(\frac{q_0^{5r-3}}{15}\Big)\!>\!\frac{q_0^{3r-3}(q_0^{2r}-1)(q_0^{3r}+1)}{c(q_0^2-1)(q_0^3+1)}=v.
\end{align*}
If $q_0=2$ and $r\geq 5$, then $c=1$ and 
\[(s+1)(s^2+1)\geq 2^{3r-3}(2^{6r-9}-2^{3r-4}+1)>\frac{2^{8r-6}}{3}>\frac{2^{3r-3}(2^{2r}-1)(2^{3r}+1)}{27}=v.\]
If $q_0=2$ and $r=3$, then $c=3$ and $(s+1)(s^2+1)\geq 30784>25536=v$. In all cases, we have shown that $(s+1)(s^2+1)>v$, which is a contradiction.

Finally, we consider the case when $q$ is even and $s\neq t$. 

Without loss of generality, we assume $s>t$. Then $T_\a=K.c$ with $K\cong\PSU(3,q_0)$ and $T_l\cong\PSU(3,q_1)$, where $q=q_0^{r_0}=q_1^{r_1}$, $r_1>r_0\geq 3$ are odd primes, and $c=\gcd(3,\frac{q+1}{q_0+1})$. 

We claim that there exists an element of order 4 contained in $T_\beta\cap T_l$ for some $\beta\in\mc{P}$.

If $T$ has a unique conjugacy class of elements of order 4, then for any $g\in T_\a$ and $h\in T_l$ with $|g|=|h|=4$, there exists $x\in T$ such that $g^x=h$. Thus $h\in T_\a^x\cap T_l=T_\beta\cap T_l$ where $\beta=\a^x$. 

If $T$ has 3 conjugacy classes of elements of order 4, then $q\equiv 2\pmod{3}$. Since $r_1$ is odd, then $q_1\equiv 2 \pmod{3}$ and so $T_l$ also has 3 conjugacy classes of elements of order 4 with representatives say $b_1,b_2$ and $b_3$, respectively. By the proof of Lemma \ref{ffqeven}, we see that $b_1,b_2$ and $b_3$ are pairwise non-conjugate in $T$. Hence for any $g\in T_\a$, there exists $x\in T$ such that $g^x=b_i$ for some $i$. Let $h=g^x$. Thus $h\in T_\a^x\cap T_l=T_\beta\cap T_l$ where $\beta=\a^x$. 

We have shown that there is an element $h\in T_\beta\cap T_l$ with $|h|=4$. By Lemma \ref{ffqeven}, we have $h^T\cap T_\beta=h^{T_\beta}$ and $h^T\cap T_l=h^{T_l}$. Recall from the proof of Lemma \ref{ffqeven} that $|C_{T_\beta}(h)|=c\cdot|C_K(h)|$. By Lemma \ref{order4}\ref{ord4.4}, we obtain $|C_T(h)|/|C_{T_\beta}(h)|=\frac{q^2}{cq_0^2}>4$ and $|C_T(h)|/|C_{T_l}(h)|=\frac{q^2}{q_1^2}>4$. It then follows from Corollary \ref{centralizer} that $C_T(h)$ is nonabelian. This contradicts the fact that $C_T(h)=\C_4^f$ given by Lemma \ref{order4}\ref{ord4.4}.

Therefore $T_\a$ cannot be isomorphic to $\PSU(3,q_0).\gcd(3,\frac{q+1}{q_0+1})$.
\end{proof}

\begin{lemma}\label{gg}
The subgroup $T_{\a}$ cannot have type (g), i.e., $T_\a\not\cong\hat~3_{+}^{1+2}\!:\!\rQ_8.\frac{\gcd(9,q+1)}{3}$, where
$q=p\equiv 2\pmod 3$ and $q\geq 11$. 
\end{lemma}

\begin{proof}
Suppose to the contrary that $T_{\a}\cong\hat~3_{+}^{1+2}\!:\!\rQ_8.\frac{\gcd(9,q+1)}{3}$. Then $T_{\a}\cong \PSU(3,2)$ if $\gcd(9,q+1)=3$, and
$T_{\a}\cong \ASL(2,3)$ if $\gcd(9,q+1)=9$.

Firstly, we assume $\gcd(9,q+1)=3$.
By \cite[Table 8.5]{lowdim}, we see that $T$ has a unique conjugacy class of subgroups isomorphic to $\PSU(3,2)$.
So we may assume $T_{\a}=T_l$. Let $g\in T_\alpha$ be of order 3. By Lemma \ref{order3}, we have $C_T(g)=Y\!:\!\C_3$ where $Y\cong \C_\frac{q+1}{3}\times \C_{q+1}$. Moreover, Magma \cite{Magma} shows that $C_{T_\alpha}(g)\cong \C_3^2$. Thus we have $\frac{|C_T(g)|}{|C_{T_\alpha}(g)|}=(q+1)^2/9$. Note that $\frac{|Y|}{|Y_\alpha|}\leq \frac{|C_T(g)|}{|C_{T_\alpha}(g)|}$, so $|Y_\alpha|\geq 3$. On the other hand, since $(|Y_\alpha|)_3\leq (|Y|)_3=3$ and $Y_\alpha\leq C_{T_\alpha}(g)\cong\C_3^2$, we conclude that $Y_\alpha\cong\C_3$. This implies that $\frac{|Y|}{|Y_\alpha|}=\frac{|C_T(g)|}{|C_{T_\alpha}(g)|}=(q+1)^2/9>4$. Finally, since $T_\alpha$ has a unique conjugacy class of elements of order 3, then $g^T\cap T_\alpha=g^{T_\alpha}$. By Corollary \ref{centralizer}, we get $Y$ is nonabelian, which is a contradiction.

Next we assume $\gcd(9,q+1)=9$. By Lemma \ref{-1square}, we have $q=p\equiv 1\pmod 4$.
In this case, Table \ref{tab_ord4} shows that $T$ has a unique conjugacy class of elements of order $4$. Take $g\in T_{\a}$ and $h\in T_l$ with $|g|=|h|=4$. There exists $x\in T$
such that $g^x=h$. Thus $h\in T_{\a}^x\cap T_l=T_{\b}\cap T_l$ where $\b=\a^x$. By Lemma \ref{order4}\ref{ord4.1}, we have $C_T(h)=\C_{(q^2-1)/3}$. Note that $T_{\b}\cong \ASL(2,3)$ has a unique conjugacy class of elements of order $4$ and $C_{T_\b}(h)\cong\C_4$. Thus $g^T\cap T_\b=g^{T_\b}$, $g^T\cap T_l=g^{T_l}$, and $|C_T(h)|/|C_{T_\b}(h)|=|C_T(h)|/|C_{T_l}(h)|=(q^2-1)/12>4$. By Corollary \ref{centralizer}, we have $C_T(h)$ is nonabelian, a contradiction.
\end{proof}

\begin{lemma}\label{hhii}
 The subgroup $T_{\a}$ cannot have type (h) or (i), i.e., $T_\a\not\cong\PSL(2,7)$ with $q\neq 5$ and $q=p\equiv 3,5,6 \pmod 7$,
and  $T_\a\not\cong\A_6$ with $q=p\equiv 11, 14 \pmod {15}$.
\end{lemma}

\begin{proof}
Suppose $T_\a\cong \PSL(2,7)$ or $\A_6$. By Lemma \ref{-1square}, we have $q\equiv1\pmod{4}$. Then by Table \ref{tab_ord4}, $T$ has a unique conjugacy class of elements of order $4$. So there exists some element $g\in T_\a\cap T_l$ of order 4 for some $\a\in\mc{P}$ and $l\in\mc{L}$. By Lemma \ref{order4}, we have $C_T(g)=\C_{(q^2-1)/d}$ where $d=\gcd(3,q+1)$. We checked by Magma \cite{Magma} that both $\PSL(2,7)$ and $\A_6$ have a unique conjugacy class of elements of order $4$, and $C_{T_\a}(g)=\C_4$. Thus $g^T\cap T_\a=g^{T_\a}$ and $|C_T(g)|/|C_{T_\a}(g)|=\frac{q^2-1}{4d}>4$. Similarly we have $g^T\cap T_l=g^{T_l}$ and $|C_T(g)|/|C_{T_l}(g)|=\frac{q^2-1}{4d}>4$. Then Corollary \ref{centralizer} implies that $C_T(g)$ is nonabelian which is a contradiction.
\end{proof}

\begin{proof}[Proof of Theorem \ref{main}]
Let $G$ be an automorphism group of the generalized quadrangle $\mc{S}$ acting primitively on both points and lines and suppose the socle of $G$ is $T=\PSU(3,q)$. Then for an arbitrary point $\a$, we have $G_\a$ is maximal in $G$ and so $T_\a=T\cap G_\a$ is isomorphic to one of the groups listed in Table \ref{tab_sub}. However, we have shown in Lemmas \ref{a}, \ref{q=5}, \ref{b}-\ref{cc}, and \ref{ff}-\ref{hhii} that $T_\a$ is not isomorphic to any of the groups in Table \ref{tab_sub}. This completes the proof.
\end{proof}

\section*{Acknowledgments} 
We would like to thank the anonymous referees whose helpful comments have improved
this paper.

\section*{Data availability} This study has no associated data.

\bibliographystyle{plain}

\begin{thebibliography}{77}

\bibitem{BE2021}
J. Bamberg, J. Evans, 
\newblock No sporadic almost simple group acts primitively on the points of a generalized quadrangle, \newblock {\it Discrete Math.} {\bf 344} (2021), 112291.

\bibitem{BGMRS2012}
J. Bamberg, M. Giudici, J. Morris, G. F. Royle, P. Spiga, 
\newblock Generalized quadrangles with a group of automorphisms acting primitively on points and lines, 
\newblock {\it J. Combin. Theory Ser. A} {\bf 119} (2012),  1479--1499.

\bibitem{BGPP2016}
J. Bamberg, S. P. Glasby, T. Popiel, C. E. Praeger, 
\newblock Generalized quadrangles and transitive pseudo-hyperovals, 
\newblock {\it J. Combin. Des.} {\bf 24} (2016), 151--164.

\bibitem{BLS2018}
J. Bamberg, C. H. Li, E. Swartz, 
\newblock A classification of finite antiflag-transitive generalized quadrangles, 
\newblock {\it Trans. Amer. Math. Soc.} {\bf 370} (2018), 1551--1601.

\bibitem{BLS2021}
J. Bamberg, C. H. Li, E. Swartz, 
\newblock A classification of finite locally 2-transitive generalized quadrangles,  
\newblock {\it Trans. Amer. Math. Soc.} {\bf 374} (2021), 1535--1578.

\bibitem{BPP2017}
J. Bamberg, T. Popiel, C. E. Praeger, 
\newblock Point-primitive, line-transitive generalized quadrangles of holomorph type, 
\newblock {\it J. Group Theory} {\bf 20} (2017), 269--287.

\bibitem{Magma}
W. Bosma, J. Cannon, C. Fieker, A. Steel, 
\newblock {\it Handbook of Magma Functions}, 2017.

\bibitem{lowdim}
J. N. Bray, D. F. Holt, and C. M. Roney-Dougal, 
\newblock {\it The maximal subgroups of the low-dimensional finite classical groups, London Mathematical Society Lecture Note Series 407}, Cambridge University  Press, Cambridge, 2013.

\bibitem{BVM1994} 
F. Buekenhout, H. Van Maldeghem,
\newblock Finite distance-transitive generalized polygons,
\newblock {\it Geom. Dedicata} {\bf 52} (1994), 41--51.


\bibitem{DM1996} 
J. D. Dixon, B. Mortimer, 
\newblock {\it Permutation groups, Grad. Texts in Math., vol. 163}, Springer-Verlag, New York, 1996.

\bibitem{FH1964}
W. Feit, G. Higman, 
\newblock The nonexistence of certain generalized polygons, 
\newblock {\it J. Algebra} {\bf 1} (1964), 114--131.


\bibitem{FL2022}
T. Feng, J. Lu,
\newblock On finite generalized quadrangles with $\mathrm{PSL}(2,q)$ as an automorphism group, 
\newblock {\it Des. Codes Cryptogr.} {\bf 91} (2023), 2347--2364. 

\bibitem{Ghi1992} 
D. Ghinelli,
\newblock Regular groups on generalized quadrangles and nonabelian difference sets with multiplier $-1$,
\newblock {\it Geom. Dedicata} {\bf 41} (1992), 165--174.

\bibitem{Kantor}
W. M. Kantor, 
\newblock Automorphism groups of some generalized quadrangles, in: {\it Advances
in Finite Geometries and Designs}, Oxford Sci. Publ., Oxford Univ. Press, New
York, 1991, pp. 251--256.

\bibitem{LS1991} 
M. W. Liebeck, J. Saxl, 
\newblock Minimal degrees of primitive permutation groups, with an application to monodromy groups of covers of Riemann surfaces, 
\newblock {\it Proc. London Math. Soc. (3)} {\bf 63} (1991), 266--314.

\bibitem{PT09}
S. E. Payne, J. A. Thas, 
\newblock {\it Finite generalized quadrangles, Second edition, EMS Series of Lectures in Mathematics}, European Mathematical Society (EMS), Z\"urich, 2009.

\bibitem{SF1973} 
W. A. Simpson, J. S. Frame, 
\newblock The character tables for $\SL(3,q)$, $\SU(3,q^2)$, $\PSL(3,q)$, $\PSU(3,q^2)$, 
\newblock {\it Canadian J. Math.} {\bf 25} (1973), 486--494.


\bibitem{Tits59} 
J. Tits, 
\newblock Sur la trialit\'{e} et certains groupes qui s'en d\'{e}duisent, 
\newblock {\it Inst. Hautes \'{E}tudes Sci. Publ. Math.} {\bf 2} (1959), 13--60.

\bibitem{Van98}
H. Van Maldeghem,
\newblock {\it Generalized polygons, Modern Birkh\"{a}user Classics},
Birkh\"{a}user/Springer Basel AG, Basel, 1998.
\end{thebibliography}

\end{document}